\documentclass[journal]{IEEEtran}
\usepackage{etex}
\usepackage[cmex10]{amsmath}
\usepackage{mathrsfs}
\usepackage{dsfont}
\usepackage{graphicx}
\usepackage{url, multicol, multirow, color, tikz, enumerate, verbatim, amsfonts}
\usepackage{threeparttable}
\usepackage{algorithmic}
\usepackage[ruled,vlined]{algorithm2e}
\usepackage{wrapfig}
\usepackage{tabularx,booktabs,caption,ragged2e}

\usepackage{tikz}
\usetikzlibrary{shapes.geometric}
\usetikzlibrary{positioning}
\usepackage{pgfplots}

\newtheorem{proposition}{Proposition}

\newcommand{\pred}{\textcolor{black}}

\newcommand{\sredVeight}{\textcolor{black}}
\newcommand{\sred}{\textcolor{black}}

\begin{document}

\title{{\color{black}Distributionally Robust Distribution Network Configuration Under Random Contingency}}
\author{Sadra Babaei,~\IEEEmembership{Student Member,~IEEE;}
Ruiwei~Jiang,~\IEEEmembership{Member,~IEEE;}
Chaoyue~Zhao,~\IEEEmembership{Member,~IEEE}
\thanks{\textcolor{black}{Sadra Babaei and Chaoyue~Zhao} are with the School of Industrial Engineering and Management, Oklahoma State University, Stillwater, OK 74074, USA. E-mail: sadra.babaei@okstate.edu\textcolor{black}{; chaoyue.zhao@okstate.edu}}
\thanks{Ruiwei~Jiang is with the Department of Industrial and Operations Engineering, University of Michigan, Ann Arbor, Michigan, 48109, USA. E-mail: ruiwei@umich.edu.}
\thanks{This work is supported in part by the National Science Foundation (CMMI-1662774, CMMI-1662589).}
}


\maketitle

\begin{abstract}
Topology design is a critical task for the reliability, economic operation, and resilience of distribution systems. This paper proposes a distributionally robust optimization (DRO) model for designing the topology of a new distribution system facing random contingencies (e.g., imposed by natural disasters). The proposed DRO model optimally configures the network topology and integrates distributed \sred{generation} to effectively meet the loads. Moreover, we take into account the uncertainty of contingency. Using the moment information of distribution line failures, we construct an ambiguity set \sred{of} the contingency probability \sred{distribution}, and minimize the expected amount of load shedding with \sred{regard to} the worst-case distribution within the ambiguity set. As compared with a classical robust optimization model, the DRO model explicitly considers the contingency uncertainty and so provides a less conservative configuration, yielding a better out-of-sample performance. We recast the proposed model to facilitate the column-and-constraint generation algorithm. We demonstrate the out-of-sample performance of the proposed approach in numerical case studies.
\end{abstract}

\section*{Nomenclature}
\addcontentsline{toc}{section}{Nomenclature}
\begin{IEEEdescription}[\IEEEusemathlabelsep\IEEEsetlabelwidth{$V_1,V_2$}]
	\item [\textcolor{black}{A. Sets}]	
	\item[$\mathcal{T}$] Set of time periods.
	\item[$\mathcal{N}$] Set of nodes.
	\item[$\mathcal{E}$] Set of power lines.
	\vspace{0.2cm}

	\item [\textcolor{black}{B. Parameters}]
	\item[$\pred{B_y}$] Available budget for power line constructions.
	\item[$\pred{B_w}$] Available number (budget) of distributed generators for allocation.
	\item[$\pred{N_z}$] Maximum number of affected power lines during the contingency.
	\item[$\pred{c_{mn}}$] Construction cost of line $(m,n)$.
	\item[$\phi_{mn}$] Resistance of the power line $(m,n)$.
	\item[$\eta_{mn}$] Reactance of the power line $(m,n)$.
	\item[$\pred{K_{mn}}$] Upper limit of active power flow in line $(m,n)$.
	\item[$\pred{R_{mn}}$] Upper limit of reactive power flow in line $(m,n)$.
	\item[$\pred{D_{nt}^p}$] Active power load at node $n$ in time $t$.
	\item[$\pred{D_{nt}^q}$] Reactive power load at node $n$ in time $t$.
	\item[$\pred{C_{n}^p}$] Active power capacity of substation or distributed generation unit at node $n$.
	\item[$\pred{C_{n}^q}$] Reactive power capacity of substation or distributed generation unit at node $n$.
	\item[$\nu^{max}$] Upper bound of voltage.
	\item[$\nu^{min}$] Lower bound of voltage.
	\item[$V_0$] Reference voltage value.
	\item[$\mu_{mn,t}^{max}$] Upper bound of failure rate in line $(m,n)$ in time $t$.
    \item[$\tau^{\mbox{\tiny rst}}_{mn}$] Minimum restoration time of line $(m,n)$ during contingency.
	\vspace{0.2cm}
	
	\item [\textcolor{black}{C. First-stage Decision Variables}]
	\item[$\pred{y_{mn}}$] Binary variable for network configuration; equals 1 if line $(m,n)$ is constructed, 0 otherwise.
	\item[$\pred{w_{n}}$] Binary variable; equals 1 if the distributed generation unit is placed at node $n$, 0 otherwise.
	\item[$\pred{f_{mn}}$] Fictitious flow across line $(m,n)$ for configuring the network.
	\item[$\mathbf{g}$] Vector of first stage decision variables including $\pred{y_{mn}}$, $\pred{w_{n}}$, and $\pred{f_{mn}}$.
	\item[$\boldsymbol{\beta}, \gamma$] Dual variables in the reformulation of the distributionally robust model.
	\vspace{0.2cm}
	
	\item [\textcolor{black}{D. Second-stage Decision Variables}]
	\item[$\pred{p_{mn,t}}$] Active power flow across line $(m,n)$ in period $t$.
	\item[$\pred{q_{mn,t}}$] Reactive power flow across line $(m,n)$ in period $t$.
	\item[$\pred{x_{nt}^p}$] Active power generation at node $n$ in period $t$.
	\item[$\pred{x_{nt}^q}$] Reactive power generation at node $n$ in period $t$.
	\item[$\nu_{nt}$] Voltage magnitude at node $n$ in period $t$.
	\item[$\pred{s_{nt}}$] Load shedding at node $n$ in period $t$.
	\item[$\boldsymbol{\pi}$] Dual variables in subproblem reformulation.
	\item[$\mathbf{u}$] Vector of second-stage decision variables including $\pred{p_{mn,t}}$, $\pred{q_{mn,t}}$, $\pred{x_{nt}^p}$, $\pred{x_{nt}^q}$, $\nu_{nt}$, and $\pred{s_{nt}}$.
	\vspace{0.2cm}
	
	\item [\textcolor{black}{E. Random Parameter}]
	\item[$\pred{z_{mn,t}}$] Bernoulli random variable; equals 0 if line $(m,n)$ is affected in period $t$, 1 otherwise.

\end{IEEEdescription}

\section{Introduction}\label{sec:intro}
Recently U.S. has witnessed repeated severe power outages due to \sred{natural disasters} such as \sred{hurricane} Sandy \cite{salman2015evaluating} and \sred{tropical storm} Irene \cite{che2014only}. Only between years of 2003--2012, nearly 679 weather-related power outages happened in the U.S. and each influenced more than \sred{50,000} customers \cite{ofeconomic}. \sred{Unfortunately}, the severity and frequency of \sred{natural disasters} have been trending \sred{upwards. For example,} in the last ten years, the U.S. has suffered \sred{from} seven of \sred{the} ten most costly storms in \sred{its} history \cite{ton2015more}.
\sred{The growing threat from natural disasters calls for better planning of the power grids to improve system resiliency.}
According to the report \cite{report} by \sred{the President's Council of Economic Advisers and the U.S. Department of Energy,} nearly 80--90\% of outages in the power system occurs along distribution systems, often leading to interruptions of power supply to end customers.
Practically, a distribution system is operated in a radial topology so as to make the design and protection coordination as simple as possible. Despite its simplicity, any contingency in the distribution system can interrupt the continuity of power supply to all customers downstream \sred{the} \sred{on-contingency area.}

The distribution network planning is widely investigated in existing \sred{literature} and broadly categorized in three \sred{parts}: distribution configuration planning (\cite{moreira2011large, kumar2014design}), distribution reconfiguration and self-healing planning (\cite{wang2015self, arefifar2013comprehensive}), and distribution reinforcement and expansion planning (\cite{yuan2016robust, ma2016resilience}).
The main objective of distribution configuration planning is to design a new system to meet the demand in the most cost-effective and reliable way.
Distribution reconfiguration and self-healing planning aim at improving or recovering of network functionality by altering the topological structure of the network. \sred{In particular, a} self-healing process is brought up when a contingency occurs in the system.
Distribution reinforcement and expansion planning involve enhancing the resilience of the network to protect against possible damages or expanding current facilities \sred{to increase reliability.}
\sred{This paper focuses} on the distribution network configuration part. \sred{As} we borrow some ideas from \sred{the} self-healing and reinforcement planning literature, we also briefly review the relevant works in these domains.

\sred{Existing} mathematical \sred{models of} distribution network configuration \sred{involve various design variables that} usually include \sred{the} location \cite{samui2012direct} and size \cite{navarro2009large} of equipments like substations \sred{and} feeders.
\sred{As the penetration of distributed generation (DG) resources grows, the location and sizing of the DG units has also received increasing attention in the literature} (see, e.g., \sred{\cite{mahmoud2016optimal}, \cite{pereira2016optimal},} \cite{atwa2010optimal, liu2011optimal}).
\sred{The network topology is another important design variable} (see, e.g., \cite{miguez2002improved, abdelaziz2010distribution, lavorato2012imposing}, \sred{\cite{el2016optimal}}). \cite{miguez2002improved} \sred{proposes an optimal network topology design that} minimizes investment and variable costs associated with power losses and reliability.
\cite{abdelaziz2010distribution} \sred{considers network reconfiguration and maintains a radial network topology by ensuring that the node-incidence matrix has non-zero determinant.} \cite{lavorato2012imposing} explicitly incorporates \sred{the} radiality constraints in the distribution \sred{system configuration} model \sred{and considers the integration of DG units.}
\sred{None} of the above works \sred{incorporate} the possibility of contingency occurrences \sred{in} the planning stage.

Most of \sred{existing} planning models in the literature \sred{incorporate contingencies in} a post-outage recovery \sred{formulation that identifies an} optimal \sred{network} reconfiguration and promptly \sred{restores} the system.
\cite{wang2015self} studies a comprehensive framework for the distribution system in both normal operation and self-healing modes. In the normal operation mode, the objective is to minimize the operation costs. When a contingency happens, the system enters the self-healing mode by sectionalizing the on--outage zone into a set of self-supplied microgrids (MGs) to \sred{pick up the maximum amount of loads.}
\cite{arefifar2013comprehensive} develops a systematic framework including planning and operating stages for a smart distribution system. In the planning stage, the goal is \sred{to construct} self-sufficient MGs using \sred{various} DGs and storage units. \sred{In} the operating stage, a new formulation that incorporates both emergency reactions and system restoration is addressed for carrying out \sred{optimal} self-healing control actions.
\cite{li2014distribution} proposes a graph-theoretic distribution system restoration algorithm to find \sred{an} optimal \sred{network reconfiguration} after \sred{multiple contingencies arise} in the system, \sred{where the} MGs are modeled as virtual feeders and the distribution system is \sred{modeled as} a spanning tree.
All of the above works are under \sred{the} premise that the contingencies have already been \sred{located} and then \sred{we perform} system reconfiguration to enhance its reliability. \sred{In contrast, this paper considers the stochasticity of the contingency (e.g., caused by natural disasters).}

\sred{Existing distribution reinforcement planning models consider stochastic contingencies and carry out} pre-event enhancement activities \sred{including} vegetation management, pole refurbishments, and undergrounding of power lines \cite{salman2015evaluating}.
\cite{yuan2016robust} presents a two-stage robust optimization \sred{model for optimally allocating DG resources and hardening lines before the upcoming natural disasters.} 
A new uncertainty set for contingency occurrences is developed to capture the spatial and temporal \sred{dynamic of hurricanes.}
\cite{ma2016resilience} \sred{proposes a new tri-level optimization approach} to mitigate \sred{the impacts} of extreme weather events on the distribution system, with the objective of minimizing hardening investment and the worst-case load shedding cost. \sred{An} infrastructure fragility model is exploited by considering \sred{a} time-varying uncertainty set of disastrous events.
\sred{Even though the above works adopt realistic uncertainty sets for modeling the contingency, challenges still exist for the robust optimization approaches. Indeed,} they completely neglect the probabilistic characteristics of \sred{the} contingency.
\sred{Accordingly, the robust optimization approaches may only focus on the worst-case contingency and yield over-conservative solutions.}

To tackle these \sred{challenges,} distributionally robust \sred{(DR)} models have been \sred{proposed} \cite{delage2010distributionally}. \sred{The DR models consider a set of probability distributions of the uncertain parameters (termed ambiguity set) using certain statistical characteristics (e.g., moments).} Then, \sred{we search for a solution that is optimal with respect to the worst-case probability distribution within the ambiguity set. DR models have been applied} on various power system problems, such as unit commitment \cite{xiong2017distributionally}, reserve scheduling \cite{bian2015distributionally}, congestion management \cite{qiu2015distributionally}, and transmission expansion planning \cite{bagheri2017data}.

To the best of our knowledge, this paper conducts the first study of DR models for distribution network configuration when facing contingency. Our main contributions include: (a) by incorporating the contingency probability distribution, our DR model is able to capture the contingencies with lower probability but high impacts, two key features of natural disaster-induced outages; (b) we recast the DR model as a two-stage robust optimization formulation that facilitate the column-and-constraint generation algorithm (see Proposition \ref{prop:ref}); (c) solving the DR model yields a worst-case contingency distribution, which can be used (e.g., in simulation models) to examine other topology configuration/re-configuration policies facing random contingency (see Proposition \ref{prop:wc-distribution}); (d) numerical case studies demonstrate the better out-of-sample performance of our DR model.


The remainder of the paper is organized as follows. In Section~\ref{sec:math}, we describe \sred{the DR formulation} including \sred{the} network configuration, \sred{the} restoration process, and \sred{the ambiguity set of} contingency \sred{probability distribution.} In Section~\ref{sec:solution}, we \sred{derive an equivalent reformulation and employ the} column-and-constraint generation framework to solve the problem. Finally, in Section~\ref{sec:case}, we \sred{conduct} case studies and \sred{analyze the} computational results.

\section{Mathematical model}\label{sec:math}
We propose a distributionally robust optimization model for \sred{a} distribution network facing random contingency. The model involves two stages. In the first stage, we form a set of \sred{radially} configured networks, each energized by a substation within the network. \sred{In addition, we allocate} a set of available DGs in the system.
Then, the contingency launches a set of disruptions to the system to inflict damages.
\sred{In} the second stage, \sred{we take} restoration actions to minimize the load shedding by rescheduling the output of substations and DGs.

\subsection{Distribution network configuration}
We plan to establish a distribution system in a new community without existing facilities. In this community, only the locations of loads and substations are identified. It is assumed that the substations are connected to a higher-level substation in the grid. Let graph $G=(\mathcal{N},\mathcal{E})$ represent the distribution network, where $\mathcal{N}$ \sred{denotes} the set of nodes and $\mathcal{E}$ denotes the set of \sred{distribution} lines that can be constructed. Also, assume that substations are located \sred{in} the set $\mathcal{R} \subset \mathcal{N}$. In the devised \sred{network configuration}, the distribution system \sred{consists of} a set of radial networks in the sense that each load bus \sred{is} connected to a substation directly or via other nodes.
\begin{figure}[!h]
	\begin{center}
		\usepgflibrary{arrows}
		\usetikzlibrary{arrows}
		\usepgflibrary{decorations.pathmorphing} 
		\usetikzlibrary{decorations.pathmorphing} 
		 {
		\begin{tikzpicture}
		[bus/.style={circle,inner sep=0pt,minimum size=4pt,fill=black},
		source/.style={circle,inner sep=0pt,minimum size=4pt,fill=red},
		source2/.style={circle,draw=purple!100,thick,
			inner sep=0pt,minimum width=2.5mm,minimum height=2.5mm},
		point/.style={circle,inner sep=0pt,minimum size=.1pt,fill=black},
		substation/.style={rectangle,draw=orange!100,thick,
			inner sep=0pt,minimum width=2.5mm,minimum height=2.5mm}]
		\node[substation] (bus1)   at (1, 1)  [label=below:\scriptsize{}] {};
		\node[bus] (bus1)   at (1, 1)  [label=below:\scriptsize{}] {};
		\node[bus] (bus2)   at (2, 1)  [label=below: \scriptsize{}] {};
		\node[substation] (bus3)   at (3, 1)  [label=below:\scriptsize{}] {};
		\node[bus] (bus3)   at (3, 1)  [label=below:\scriptsize{}] {};
		\node[bus] (bus4)   at (3, 0)  [label={[label distance=-.3mm]-45:\scriptsize{}}] {};
		\node[bus] (bus5)   at (1, 0)  [label=above:\scriptsize{}] {};
		
		\node[source2] (bus6)   at (5.5, 1.5)  [label=above:\scriptsize{}] {};
		\node[source] (bus6)   at (5.5, 1.5)  [label=above:\scriptsize{$s$}] {};
		\node[substation] (bus7)   at (4.5, 1)  [label=above:\scriptsize{}] {};
		\node[bus] (bus7)   at (4.5, 1)  [label=above:\scriptsize{}] {};
		\node[bus] (bus8)   at (5.5, 1)  [label=above:\scriptsize{}] {};
		\node[substation] (bus9)   at (6.5, 1)  [label=above:\scriptsize{}] {};
		\node[bus] (bus9)   at (6.5, 1)  [label=above:\scriptsize{}] {};
		\node[bus] (bus10)   at (6.5, 0)  [label=above:\scriptsize{}] {};
		\node[bus] (bus11)   at (4.5, 0)  [label=above:\scriptsize{}] {};
		
		\node[substation] (bus12)   at (1.3, -.55)  [label=below:\scriptsize{$\mbox{Substation}$}] {};
		\node[source2] (bus13)   at (3.2, -.6)  [label=below:\scriptsize{}] {};
		\node[source] (bus13)   at (3.2, -.6)  [label=below:\scriptsize{$\mbox{\sred{Higher-level} substation}$}] {};
		\node[point] (bus14)   at (4.8, -.65)  [label=below:\scriptsize{}] {};
		\node[point] (bus15)   at (5.4, -.65)  [label=below:\scriptsize{$\mbox{New lines}$}] {};
		\draw [dashed] (bus14) -- (bus15);
		\draw [-] (bus2) -- (bus3);
		\draw [-] (bus2) -- (bus4);
		\draw [-] (bus1) -- (bus5);
		
		\draw [dashed] (bus6) -- (bus7);
		\draw [dashed] (bus6) -- (bus9);
		\draw [-] (bus8) -- (bus9);
		\draw [-] (bus8) -- (bus10);
		\draw [-] (bus7) -- (bus11);

		\path[draw=blue,solid,line width=1.3mm,fill=blue,
		preaction={-triangle 90,thin,draw=blue,fill=blue,shorten >=-1mm}
		] (3.5, .75) -- (4, .75);
		\path[draw=blue,solid,line width=1.3mm,fill=blue,
		preaction={-triangle 90,thin,draw=blue,fill=blue,shorten >=-1mm}
		] (4, .75) -- (3.5, .75);
		\end{tikzpicture}
	}
	\end{center}
	\caption{\sred{Example of a spanning tree representation}} \label{example}
\end{figure}
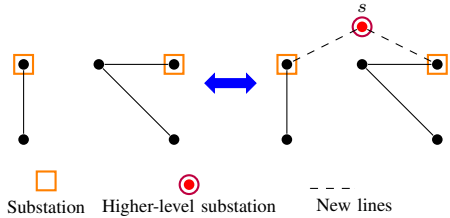
In other words, we \sred{construct} a spanning forest with $|\mathcal{R}|$ components, \sred{each rooted at} one substation. For this purpose, we add a new \sred{higher-level node $s$} to graph $G$ and connect it to all \sred{substation} nodes, i.e., nodes in $\mathcal{R}$. We call the new graph $G'=(\mathcal{N}',\mathcal{E'})$. Now constructing a \sred{spanning forest rooted} in $\mathcal{R}$ is equivalent to \sred{constructing} a spanning tree of this new graph $G'$, where all \sred{newly} added lines (i.e., $\mathcal{E'} \backslash \mathcal{E}$) are included in the tree \sred{(see Fig. \ref{example} for an example).}
To formulate the spanning tree, we employ the single commodity formulation \cite{magnanti1995optimal} as follows:
\begin{subequations}
\begin{eqnarray}
&& \hspace{-.3in}\sum_{ n| (s,n) \in \mathcal{E'}} f_{sn} = |\mathcal{N}'| -1, \label{tree1}\\
&& \hspace{-.3in}\sum_{ m| (m,n) \in \mathcal{E'}} f_{mn} - \sum_{m| (n,m) \in \mathcal{E'}} f_{nm} = 1, \ \ \forall{n \in \mathcal{N}' \backslash s},\label{tree2}\\
&& \hspace{-.3in}\sum_{(m,n) \in \mathcal{E'}} y_{mn} = |\mathcal{N}'| - 1,\label{tree3}\\
&& \hspace{-.3in}f_{mn} \leq (|\mathcal{N}'|-1) y_{mn},\ \ \forall{(m,n) \in \mathcal{E'}},\label{tree4}\\
&& \hspace{-.3in} y_{mn}=1, \ \ \forall{(m,n)\in \mathcal{E'} \backslash \mathcal{E}},\label{tree6}\\
&& \hspace{-.3in} f_{mn} \geq 0, \ \ y_{mn} \in \{0,1\}, \ \ \forall{(m,n) \in \mathcal{E'}}.  \label{tree7}
\end{eqnarray}
We remark that $f_{mn}$ does not represent the power flow along the line $(m,n)$. Instead, it represents fictitious flow to mathematically guarantee that the distribution network is radial. Constraint \eqref{tree1} indicates that there must be $|\mathcal{N}^{'}|-1$ arcs leaving the root node $s$ in order to form a spanning tree. Constraints \eqref{tree2} ensure the connectivity of the spanning tree. Constraint \eqref{tree3} specifies \sred{that,} in \sred{the constructed} spanning tree, the number of connected lines should be one unit less than the number of nodes. Constraints \eqref{tree4} designate that the capacity of fictitious flow on each line should be no more than the total number of connected lines. Constraints \eqref{tree6} indicate that all substations should be connected to the \sred{higher-level node} $s$.

Furthermore, we consider the budget constraints on the number of available DG units for installation and the total construction costs, as stated in \eqref{budgetw} and \eqref{budgety}, respectively:
\begin{eqnarray}
&& \hspace{-.0in} \sum_{n \in \mathcal{N}} w_n \leq B_w, \label{budgetw} \\
&& \hspace{-.0in} \sum_{(m,n) \in \mathcal{E}} c_{mn}y_{mn} \leq B_y. \label{budgety}
\end{eqnarray}
\end{subequations}

\subsection{Post-contingency restoration process}
In this study, we model the whole restoration process using a number of corrective actions to minimize the load shedding.
\sred{We adopt the well-studied} linearized approximation of \sred{the} DistFlow model \sred{(see, e.g., \cite{yuan2016robust, baran1989network})} to formulate power flow in \sred{the} distribution system \sred{after the contingency.} According to this model, active and reactive power balance flows at each bus are \sred{expressed} as follows:
\begin{subequations}
\begin{eqnarray}
&&\hspace{-0.0in} \sum_{k|(n,k) \in \mathcal{E}} p_{nk,t} = p_{mn,t} - D_{nt}^{p} + x_{nt}^p + s_{nt}, \nonumber\\
&&\hspace{-0.0in}\forall{n \in \mathcal{N}},\ \ \sred{\forall{(m,n)} \in \mathcal{E}}, \ \ \forall{t \in \mathcal{T}}, \label{BalanceActive}\\
&&\hspace{-0.0in} \sum_{k|(n,k) \in \mathcal{E}} q_{nk,t} = q_{mn,t} - D_{nt}^{q} + x_{nt}^q, \nonumber\\
&&\hspace{-0.0in}\forall{n \in \mathcal{N}},\ \ \sred{\forall{(m,n)} \in \mathcal{E}}, \ \ \forall{t \in \mathcal{T}}. \label{BalanceReactive}
\end{eqnarray}

According to the linearized DistFlow model, the relationship of voltage level \sred{between any pair of adjacent nodes} is characterized by the following constraints:
\begin{eqnarray}
&&\hspace{-0.0in} \nu_{nt} y_{mn} = \nu_{mt} y_{mn} - (\phi_{mn} p_{mn,t} + \eta_{mn} q_{mn,t}) / V_0, \nonumber\\
&&\hspace{-0.0in}\forall{m, n \in \mathcal{N} | (m,n) \in \mathcal{E}}, \ \ \forall{t \in \mathcal{T}}. \label{BalanceVoltage}
\end{eqnarray}

Moreover, the voltage level at each node should be within a permissible range:
\begin{eqnarray}
&&\hspace{-0.0in} \nu^{min} \leq \nu_{nt} \leq \nu^{max}, \ \ \forall{n \in \mathcal{N}}, \ \ \forall{t \in \mathcal{T}}. \label{Voltage}
\end{eqnarray}

Additionally, if line $(m,n)$ is not constructed in the configuration stage or constructed but disrupted during the contingency, the power flow on line $(m,n)$ should be zero. These restrictions are described by the following constraints:
\begin{eqnarray}
&&\hspace{-0.3in} 0 \leq p_{mn,t} \leq K_{mn} z_{mn,t}  y_{mn}, \
\sred{\forall{(m,n) \in \mathcal{E}},} \ \forall{t \in \mathcal{T}}, \label{UpperActive} \\
&&\hspace{-0.3in} 0 \leq q_{mn,t} \leq R_{mn}z_{mn,t} y_{mn}, \
\sred{\forall{(m,n) \in \mathcal{E}},} \ \forall{t \in \mathcal{T}}. \label{UpperReactive}
\end{eqnarray}

In our proposed framework, each radial network is rooted at a node where the substation is placed. Moreover, \sred{the DG units can supply power not only to their neighboring loads but also to all nodes in the connected network.}
The active and reactive power capacity of the substations and DGs are described by the following constraints:
\begin{eqnarray}
&&\hspace{-0.0in} 0 \leq x_{nt}^{p} \leq  C_n^{p}, \ \ \forall{n \in \mathcal{R}}, \ \ \forall{t \in \mathcal{T}}, \label{subCapP}\\
&&\hspace{-0.0in} 0 \leq x_{nt}^{q} \leq C_{n}^{q}, \ \ \forall{n \in \mathcal{R} }, \ \ \forall{t \in \mathcal{T}}, \label{subCapQ}\\
&&\hspace{-0.0in} 0 \leq x_{nt}^{p} \leq w_n C_n^{p}, \ \ \forall{n \in \mathcal{N} \backslash \mathcal{R}}, \ \ \forall{t \in \mathcal{T}}, \label{UpperActivegeneration}\\
&&\hspace{-0.0in} 0 \leq x_{nt}^{q} \leq w_n C_{n}^{q}, \ \ \forall{n \in \mathcal{N} \backslash \mathcal{R}}, \ \ \forall{t \in \mathcal{T}}. \label{UpperReactivegeneration}
\end{eqnarray}

Finally, the unsatisfied active demand at each node should be no more than the active demand at that node:
\begin{eqnarray}
&&\hspace{-0.0in} 0 \leq s_{nt} \leq D_{nt}^{p}, \ \ \forall{n \in \mathcal{N}}, \ \ \forall{t \in \mathcal{T}}. \label{shed}
\end{eqnarray}
\end{subequations}

\subsection{\sred{Ambiguity set of contingency}}\label{sec:ambiguity}
Different approaches have been proposed in the literature to deal with the uncertainty of \sred{contingency.}
Stochastic programming \sred{(SP) is well-known for modeling contingency due to natural disasters} (see, e.g., \cite{gao2017resilience, arab2015stochastic, yamangil2014designing}). Using statistical methods, \sred{SP estimates} the joint probability distribution of \sred{contingency} and then generates \sred{a set of scenarios to represent} the stochastic \sred{contingency in decision making.} The major drawback of this approach is that the underlying probability distribution often cannot be estimated accurately, and the computational effort \sred{significantly increases} as the number of contingency scenarios increases.
Robust optimization \sred{(RO)} is another well-known approach to cope with the \sred{uncertainty of contingency} (see, e.g., \cite{yuan2016robust, ma2016resilience}). \sred{Applied on} the distribution network configuration problem, \sred{RO identifies} the most critical contingencies by solving the following bilevel model:
\begin{subequations}
\begin{eqnarray}
&& \hspace{-.6in} \max_{\mathbf{z} \in \mathfrak{D}{(\mathbf{g})}} Q(\mathbf{g},\mathbf{z}) \label{dist0}\\
&& \hspace{-.6in} \mbox{s.t.} \ \ \mathfrak{D}{(\mathbf{g})}=\Bigg\{\sum_{(m,n) \in \mathcal{E}} (1-z_{mn,t}) \leq N_z, \ \forall{t \in \mathcal{T}}, \label{dist4}\\
&& \hspace{-.1in} 1-z_{mn,t} \leq y_{mn}, \ \forall{(m,n) \in \mathcal{E}}, \ \forall t \in \mathcal{T}, \label{dist3} \\
&& \hspace{-.1in} z_{mn,t+\tau} \leq z_{mn,t}, \ \forall{(m,n) \in \mathcal{E}}, \ \forall \tau \leq \tau^{\mbox{\tiny rst}}_{mn} \label{dist5}
\Bigg\},
\end{eqnarray}
\end{subequations}
where,
\begin{subequations}
\begin{eqnarray}
&&\hspace{-0.5in} Q(\mathbf{g},\mathbf{z})= \min_{\mathbf{u} \in \mathcal{H}(\mathbf{g},\mathbf{z})}  \sum_{t \in \mathcal{T}} \sum_{n \in \mathcal{N}} s_{nt} \label{MinShedding}, \\
&&\hspace{-0.5in} \mbox{s.t.}\ \ \mathcal{H}(\mathbf{g},\mathbf{z})= \bigg\{ \mathbf{u}:\mbox{Constraints \eqref{BalanceActive}-\eqref{shed} } \bigg\} \label{Hconstraint},
\end{eqnarray}
\end{subequations}
{\color{black}$\mathbf{g}:=(\mathbf{y}, \mathbf{w}, \mathbf{f})$ indicates network configuration and DG allocation decision variables,} \sred{$\mathbf{u}:=(\mathbf{p}, \mathbf{q}, \mathbf{x}^p, \mathbf{x}^q, \boldsymbol{\nu}, \mathbf{s})$ denotes the post-contingency decision variables, and $Q(\mathbf{g},\mathbf{z})$ represents the minimum load shedding for given topology $\mathbf{g}$ and contingency $\mathbf{z}$. Moreover, $\mathfrak{D}{(\mathbf{g})}$ specifies the set of all possible contingency scenarios.}
We assume \sred{by constraints} \eqref{dist4} that the number of simultaneous line outages is bounded by $N_z$, which can be calibrated based on reliability analyses of distribution lines during contingency (see, e.g.,~\cite{sa2002reliability}). \sredVeight{Constraints \eqref{dist3} designate that only constructed lines can be affected, i.e., $z_{mn,t}$ is set to be one whenever $y_{mn}$ equals zero. However, variables $z_{mn,t}$ only appear in constraints \eqref{UpperActive}--\eqref{UpperReactive}, whose right-hand sides equal zero if $y_{mn}=0$, regardless of the value of $z_{mn,t}$. Hence, we can relax constraints \eqref{dist3} without loss of optimality.}
Constraints \eqref{dist5} model the minimum restoration time of failing distribution lines. As discussed in Section \ref{sec:intro}, the RO model may only focus on the worst-case contingency (i.e., $\mathbf{z} \in \mathfrak{D}{(\mathbf{g})}$ that maximizes $Q(\mathbf{g},\mathbf{z})$ in \eqref{dist0}) and yield over-conservative topology design and/or DG allocation.

To overcome the challenges of \sred{the classical} stochastic and robust approaches, we propose a \sred{DR} framework \sred{considering a family of joint probability distributions of contingency based on the moment information of the random parameters (see, e.g., \cite{delage2010distributionally, xiong2017distributionally, zhao2017distributionally}).} \sred{More specifically, we define the ambiguity set as follows:
\begin{eqnarray}
&& \hspace{-.5in} \mathscr{D} = \Bigg\{\mathbb{P} \in \mathcal{P}\Big(\mathfrak{D}{(\mathbf{g})}\Big):  0 \leq  E_{\mathbb{P}} [1-{\mathbf{z}}] \leq \mathbf{\boldsymbol\mu}^{max}  \label{dist2}
\Bigg\},
\end{eqnarray}
where $\mathcal{P}\Big(\mathfrak{D}{(\mathbf{g})}\Big)$ consists of all probability distributions on a sigma-field of $\mathfrak{D}{(\mathbf{g})}$.
Constraints \eqref{dist2} imply that the marginal probability of each line $(m,n)$ not working during time unit $t$ has an upper limit $\mathbf{\boldsymbol\mu}^{max}$. We note that, although $\mathscr{D}$ models the contingency of new distribution lines, the distributional information (e.g., $\mathbf{\boldsymbol\mu}^{max}$ and $\mathfrak{D}{(\mathbf{g})}$) can be calibrated based on reliability analyses of distribution lines (see, e.g.,~\cite{sa2002reliability}). Accordingly, we consider the following DR model:}
\begin{eqnarray}
&& \hspace{-.1in} \max_{\mathbb{P} \in \mathscr{D}} E_{\mathbb{P}} [Q(\mathbf{g},{\mathbf{z}})]. \label{worstcaseproblem}
\end{eqnarray}
Here, instead of considering the worst-case scenario of contingency as \sred{in the RO model,} we consider the worst-case distribution of contingency and the corresponding expected load shedding. Hence, our approach, though \sred{still risk-averse,} is less conservative than the \sred{RO} approach.

\subsection{Distributionally robust optimization model}
Our distributionally robust optimization model aims to find \sred{an} optimal distribution system configuration to minimize the load shedding under random contingency:
\begin{subequations}
\begin{eqnarray}
&&\hspace{-0.5in} \min_{\mathbf{g} \in \mathcal{G}}  \max_{\mathbb{P} \in    \mathscr{D}} E_{\mathbb{P}}[Q(\mathbf{g},{\mathbf{z}})], \label{objective} \\
&&\hspace{-0.5in} \mbox{s.t.} \ \ \mathcal{G}= \bigg\{\mathbf{g}: \mbox{Constraints \eqref{tree1}-\eqref{budgety}} \bigg\}. \label{DRcons}
\end{eqnarray}
\end{subequations}

In above formulation, the objective function \eqref{objective} aims to minimize the worst-case expected load shedding $Q(\mathbf{g},\mathbf{z})$.
\section{Solution Methodology}\label{sec:solution}
In this section, we \sred{first derive reformulations of the worst-case expectation model \eqref{worstcaseproblem} and the DR model \eqref{objective}--\eqref{DRcons}, respectively.
Then, we describe a solution approach based on the column-and-constraint generation (CCG) \sred{framework.} Finally, we derive the worst-case distribution of contingency. }
\subsection{Problem reformulation}
\begin{proposition} \label{prop:ref}
For \sred{fixed $\mathbf{g} \in \mathcal{G}$, we have}
\begin{eqnarray}
&&\hspace{-0.5in} \max_{\mathbb{P} \in \mathscr{D}} E_{\mathbb{P}}[Q(\mathbf{g},{\mathbf{z}})]= \min_{\mathbf{\boldsymbol\beta}  \geq 0} \max_{\mathbf{z} \in \mathfrak{D}{(\mathbf{g})}} \biggl\{Q(\mathbf{g},\mathbf{z})\nonumber\\
&&\hspace{-0.5in}+\sum_{t \in \mathcal{T}} \sum_{(m,n) \in \mathcal{E}} (\mu_{mn,t}^{max} +z_{mn,t}-1) \beta_{mn,t}\biggr\}, \nonumber
\end{eqnarray}
where dual variables $\mathbf{\boldsymbol\beta}$ are associated with constraints \eqref{dist2}.
\end{proposition}

The proof is given in \sred{the} appendix.
\sred{By Proposition \ref{prop:ref} and combining two minimizations, we obtain the following equivalent reformulation of formulation \eqref{objective}--\eqref{DRcons}:}\sred{
\begin{eqnarray}
&&\hspace{-0.8in} \min_{\mathbf{\boldsymbol\beta} \geq 0, \mathbf{g} \in \mathcal{G}} \hspace{0.06in} \max_{\mathbf{z} \in \mathfrak{D}{({\mathbf{g}})}} \min_{\mathbf{u} \in \mathcal{H}({\mathbf{g}},\mathbf{z})}\ \sum_{t \in \mathcal{T}} \sum_{n \in \mathcal{N}} s_{nt}  \nonumber\\
&&\hspace{-0.3in} + \sum_{t \in \mathcal{T}} \sum_{(m,n) \in \mathcal{E}} (\mu_{mn,t}^{max} +z_{mn,t}-1) \beta_{mn,t}. \label{twostageProblem}
\end{eqnarray}
}Therefore, the DR model \eqref{objective}--\eqref{DRcons} is transformed into the classical robust optimization problem \eqref{twostageProblem}.
\subsection{\sred{Column-and-constraint generation framework}}
We employ \sred{the CCG framework} \cite{zeng2013solving} to solve \sred{the} problem \eqref{twostageProblem}. \sred{We describe the master problem in the $r^{th}$ iteration of the CCG framework as follows:}
\begin{subequations}
\begin{eqnarray}
&&\hspace{-0.5in} \min_{\mathbf{\boldsymbol\beta}\geq 0, \mathbf{g} \in \mathcal{G}, \lambda, \mathbf{u}^j}   \sum_{t \in \mathcal{T}} \sum_{(m,n) \in \mathcal{E}} (\mu_{mn,t}^{max}-1) \beta_{mn,t} + \lambda \label{master1}\\
&&\hspace{-0.5in} \mbox{s.t.} \ \ \lambda \geq \sum_{t \in \mathcal{T}} \sum_{n \in \mathcal{N}} s_{nt}^{j} + \sum_{t \in \mathcal{T}} \sum_{(m,n) \in \mathcal{E}} z_{mn,t}^{j} \beta_{mn,t}, \nonumber\\
&&\hspace{-0.2in}  \forall{\mathbf{z}^j \in \mathcal{F}}, \ \ \forall{j=1,...,r}, \label{master2}\\
&&\hspace{-0.2in} \mathbf{u}^j \in \mathcal{H}(\mathbf{g},\mathbf{z}^{j}), \ \ \forall{\mathbf{z}^j \in \mathcal{F}}, \ \ \forall{j=1,...,r}, \label{master3}
\end{eqnarray}
\end{subequations}
where $\mathcal{F} \subseteq \mathfrak{D{(\mathbf{g})}}$. \sred{In} the CCG framework, \sred{set $\mathcal{F}$ is iteratively augmented by incorporating more scenarios. Note that,} the master problem is a relaxation of the original problem, in which the set of contingency $\mathfrak{D{(\mathbf{g})}}$ \sred{consists of all possible scenarios satisfying constraints \eqref{dist4} (note that, as discussed in Section \ref{sec:ambiguity}, we have relaxed \eqref{dist3} without loss of optimality).}
Therefore, solving the master problem \eqref{master1}--\eqref{master3} yields a lower bound for that optimal value of \eqref{twostageProblem}. \sred{In contrast, the following subproblem yields an upper bound:}
\begin{eqnarray}
&&\hspace{-0.5in} \max_{\mathbf{z} \in \mathfrak{D}{(\hat{\mathbf{g}})}} \min_{\mathbf{u} \in \mathcal{H}(\hat{\mathbf{g}},\mathbf{z})} \sum_{t \in \mathcal{T}} \sum_{n \in \mathcal{N}} s_{nt} + \sum_{t \in \mathcal{T}} \sum_{(m,n) \in \mathcal{E}} \hat{\beta}_{mn,t} z_{mn,t},\label{subproblem}
\end{eqnarray}
\sred{where decisions $\hat{\mathbf{g}}$ and $\hat{\mathbf{\boldsymbol\beta}}$ are obtained from solving the master problem \eqref{master1}--\eqref{master3}. Note that ($\hat{\mathbf{g}}$, $\hat{\mathbf{\boldsymbol\beta}}$) is feasible to the problem \eqref{twostageProblem}.} Hence, the optimal objective value of \sred{ \eqref{subproblem}, plus constant $\sum_{t \in \mathcal{T}} \sum_{(m,n) \in \mathcal{E}} (\mu_{mn,t}^{max} - 1) \hat{\beta}_{mn,t}$, is} an upper bound for \eqref{twostageProblem}. Moreover, since the inner minimization \sred{problem of \eqref{subproblem}} is always feasible and bounded (a trivial solution is when all loads are shed), \sred{we take the dual of this minimization problem with strong duality} and convert the bilevel subproblem \eqref{subproblem} into the following \sred{single-level} bilinear maximization problem:
\begin{subequations}
\begin{eqnarray}
&&\hspace{-0.1in} \max_{\mathbf{z} \in \mathfrak{D}{(\mathbf{g})},\mathbf{\boldsymbol\pi}}  \sum_{t \in \mathcal{T}}
\sum_{(m,n) \in \mathcal{E}} \beta_{mn,t} z_{mn,t}
-\sum_{t \in \mathcal{T}} \sum_{n \in \mathcal{N}} D_{nt}^{p} \pi_{nt}^1  \nonumber\\
&&\hspace{-0.1in} - \sum_{t \in \mathcal{T}} \sum_{n \in \mathcal{N}} D_{nt}^{q} \pi_{nt}^2  + \sum_{t \in \mathcal{T}} \sum_{(m,n) \in \mathcal{E}} K_{mn} \pi_{mn,t}^3 z_{mn,t} y_{mn}\nonumber\\
&&\hspace{-0.1in} + \sum_{t \in \mathcal{T}} \sum_{(m,n) \in \mathcal{E}} R_{mn} \pi_{mn,t}^4 z_{mn,t}y_{mn} \nonumber\\
&&\hspace{-0.1in} +\sum_{t \in \mathcal{T}} \sum_{n \in \mathcal{N}\backslash \mathcal{R}} w_n C_n^p \pi_{nt}^5 + \sum_{t \in \mathcal{T}} \sum_{n \in \mathcal{N}\backslash \mathcal{R}} w_n C_n^q \pi_{nt}^6  \nonumber\\
&&\hspace{-0.1in} +\sum_{t \in \mathcal{T}} \sum_{n \in \mathcal{R}}  C_n^p \pi_{nt}^7 + \sum_{t \in \mathcal{T}} \sum_{n \in \mathcal{R}}  C_n^q \pi_{nt}^8  \nonumber\\
&&\hspace{-0.1in} +  \sum_{t \in \mathcal{T}} \sum_{n \in \mathcal{N}} v^{max} \pi_{nt}^{10} - \sum_{t \in \mathcal{T}} \sum_{n \in \mathcal{N}} v^{min} \pi_{nt}^{11}  \nonumber\\
&&\hspace{-0.1in} +\sum_{t \in \mathcal{T}} \sum_{n \in \mathcal{N}} D_{nt}^p \pi_{nt}^{12}  \label{sub1} \\
&& \mbox{s.t.} \ \pi_{mn,t}^3 + \pi_{mt}^1 - \pi_{nt}^1 + \frac{\phi_{mn}}{V_0}\pi_{nt}^{9} \leq 0, \nonumber\\
&&\hspace{0.2in}\forall{m, n \in \mathcal{N} | (m,n) \in \mathcal{E}}, \ \ \forall{t \in \mathcal{T}},\label{sub2}\\
&&\hspace{0.2in} \pi_{mn,t}^4 + \pi_{mt}^2 - \pi_{nt}^2 + \frac{\eta_{mn}}{V_0}\pi_{nt}^{9} \leq 0, \nonumber\\
&&\hspace{0.2in}\forall{m, n \in \mathcal{N} | (m,n) \in \mathcal{E}}, \ \ \forall{t \in \mathcal{T}},\label{sub3}\\
&&\hspace{0.2in} -\pi_{nt}^1 + \pi_{nt}^5 \leq 0, \ \ \forall{n \in \mathcal{N} \backslash \mathcal{R}}, \ \ \forall{t \in \mathcal{T}},\label{sub4}\\
&&\hspace{0.2in} -\pi_{nt}^2 + \pi_{nt}^6 \leq 0, \ \ \forall{n \in \mathcal{N}\backslash \mathcal{R}}, \ \ \forall{t \in \mathcal{T}},\label{sub5}\\
&&\hspace{0.2in} -\pi_{nt}^1 + \pi_{nt}^7 \leq 0, \ \ \forall{n \in \mathcal{R}}, \ \ \forall{t \in \mathcal{T}},\label{sub6}\\
&&\hspace{0.2in} -\pi_{nt}^2 + \pi_{nt}^8 \leq 0, \ \ \forall{n \in \mathcal{R}}, \ \ \forall{t \in \mathcal{T}},\label{sub7}\\
&&\hspace{0.2in} \pi_{nt}^{10} - \pi_{nt}^{11} + \pi_{jt}^{9} - \sum_{i| (n,i) \in \mathcal{E}}\pi_{it}^{9} \leq 0, \nonumber \\
&&\hspace{0.2in} j| (j,n) \in \mathcal{E}, \ \ \forall{n \in \mathcal{N}}, \ \ \forall{t \in \mathcal{T}},\label{sub8}\\
&&\hspace{0.2in} -\pi_{nt}^1 + \pi_{nt}^{12} \leq 1, \ \ \forall{n \in \mathcal{N}}, \ \ \forall{t \in \mathcal{T}},\label{sub9}\\
&&\hspace{-0.1in} \pi_{nt}^1, \pi_{nt}^2, \pi_{nt}^{9} \mbox{  are free and other variables are nonpositive,} \nonumber
\end{eqnarray}
\end{subequations}
where $\mathbf{\boldsymbol\pi}$ \sred{represent} dual variables \sred{pertaining} to constraints \sred{\eqref{BalanceActive}--\eqref{shed}.}
Note that bilinear terms $\mathbf{\boldsymbol\pi} \mathbf{z}$ in the objective function \sred{\eqref{sub1} }can be linearized using \sred{the} McCormick method \cite{mccormick1976computability}, \sred{which recasts the problem \eqref{sub1}--\eqref{sub9} as a mixed-integer linear program and facilitates efficient off-the-shelf solvers like CPLEX. }
The CCG framework is summarized as follows:\\
\textbf{Step 0}: \sred{Initialization.
Pick an optimality gap $\epsilon$.
Set $\mbox{LB}=-\infty$, $\mbox{UB}=+\infty$, set of contingencies $\mathcal{F}= \O$, and iteration index $r=1$.}\\
\textbf{Step 1}: Solve the master problem \eqref{master1}--\eqref{master3}, obtain \sred{the} optimal value $\mbox{objMP}$ and optimal configuration decisions $\hat{\mathbf{g}}^r$ and $\hat{\mathbf{\boldsymbol\beta}}^r$, and update $\mbox{LB}=\mbox{objMP}$.\\
\textbf{Step 2}: Solve the subproblem \eqref{sub1}--\eqref{sub9}, obtain \sred{the} optimal value $\mbox{objSP}$ and \sred{an} optimal contingency scenario $\hat{\mathbf{z}}^r$. Update $\mbox{UB}=\min\{\mbox{UB},\mbox{objSP}+\sum_{t \in \mathcal{T}}\sum_{(m,n) \in \mathcal{E}} (\mu_{mn,t}^{max}-1) \hat{\beta}_{mn,t}\}$, and $\sred{\mathcal{F} = \mathcal{F} \cup \{\hat{\mathbf{z}}^r\}}$. \\
\textbf{Step 3}: If $\mbox{Gap}=(\mbox{UB}-\mbox{LB})/\mbox{LB} \leq \epsilon$, then terminate and output $\hat{\mathbf{g}}^r$ as an optimal solution; otherwise, update $r=r+1$ and go to the next step.\\
\textbf{Step 4}: Create second-stage variables $\mathbf{u}^r$ and the corresponding constraints $\mathbf{u}^r \in \mathcal{H}(\mathbf{g},\hat{\mathbf{z}}^r)$. \sred{Add} them to the master problem and go to \sred{Step} 1.\\

An important by-product of the CCG framework is the worst-case contingency probability distribution, which is formalized in the following proposition. The proof is given in \sred{the} appendix.
\begin{proposition}\label{prop:wc-distribution}
Suppose that the CCG framework terminates at the $R^{\mbox{\tiny th}}$ iteration with optimal solutions $(\hat{\beta}^R, \hat{g}^R, \hat{\lambda}^R, \{\hat{u}^j\}_{j = 1, \ldots, R})$. Then, if we resolve formulation \eqref{master1}--\eqref{master3} with variables $g$ and $u^j$ fixed at $\hat{g}^R$ and $\hat{u}^j$, respectively, then the dual optimal solutions associated with constraints \eqref{master2}, denoted as $\{\psi_j\}_{j = 1, \ldots, R}$, characterize the worst-case contingency probability distribution, i.e., $\mathbb{P}\{\mathbf{z} = \mathbf{z}^j\} = \psi_j$, $\forall j = 1, \ldots, R$.
\end{proposition}


\section{Case Study}\label{sec:case}
To evaluate the effectiveness of our approach, we \sred{conduct} two \sred{case studies}. In the first \sred{study}, the distribution network includes 33 nodes, 3 substations, and 2 DG units for \sred{allocation}. The 3 substations are located at nodes 1, 11, and 25, \sred{respectively}. In the second \sred{study}, the system contains 69 nodes, 4 substations, and 3 DG units for \sred{allocation}. The substations are located at nodes 1, 13, 39, and 61, \sred{respectively.} The active and reactive power \sred{capacities} of the DGs are assumed to be $100 \sred{\mbox{KW}}$ and $50 \sred{\mbox{KVar}}$, respectively. In both studies, we consider 24 hours in the post-contingency restoration, i.e., $\mathcal{T} = \{1, \ldots, 24\}$. The active and reactive power loads at each node are randomly generated from intervals $[30,200] \sred{\mbox{KW}}$ and $[5,100] \sred{\mbox{KVars}}$, respectively. The construction costs for \sred{distribution} lines \sred{are randomly generated} from intervals proportional to their length. \sred{Overall, the construction costs are within the interval} \sred{$\$[40,100]\times10^4$.} The contingency status for \sred{distribution} lines is assumed to follow independent Bernoulli distributions with different failure \sred{probabilities} that vary within the interval $[0,0.01]$. Unless stated otherwise, we set the construction budget $B_y$ and the maximum number of affected lines $N_z$ to be \sred{$\$1770\times10^4$} and $3$ \sred{for the 33-node system, and $\$4480\times10^4$ and $4$ for the 69-node system,} respectively. \sred{All case studies} are implemented in C++ with \sred{CPLEX} 12.6 on a computer with Intel Xeon 3.2 GHz and 8 GB memory.

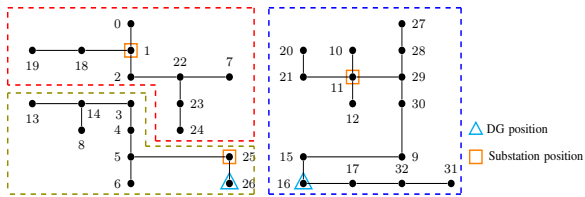
\begin{figure}[!h]
	\begin{center}
		\usepgflibrary{arrows}
		\usetikzlibrary{arrows}
		\usepgflibrary{decorations.pathmorphing} 
		\usetikzlibrary{decorations.pathmorphing} 
		\resizebox {7.9cm}{2.5cm} {
		\begin{tikzpicture}
		[bus/.style={circle,inner sep=0pt,minimum size=4pt,fill=black},
		point/.style={circle,inner sep=0pt,minimum size=2pt,fill=black},
		substation/.style={rectangle,draw=orange!100,thick,
			inner sep=0pt,minimum width=2.5mm,minimum height=2.5mm},
		load/.style={->,shorten >=0.5pt,>=triangle 45,semithick},
		dg/.style={regular polygon, regular polygon sides=3,inner sep=2pt ,draw=cyan!100,fill=white, thick},legend1/.style={regular polygon, regular polygon sides=3,inner sep=1.5pt ,draw=cyan!100,fill=white, thick},legend2/.style={rectangle,draw=orange!100,fill=white!100,thick,
			inner sep=0pt,minimum width=2.5mm,minimum height=2.5mm}]
		\node[bus] (bus17)   at (6.5, 6.5)  [label=above:\scriptsize{$17$}] {};
		\node[bus] (bus32)   at (7.5, 6.5)  [label=above: \scriptsize{$32$}] {};
		\node[bus] (bus31)   at (8.5, 6.5)  [label=above:\scriptsize{$31$}] {};
		\node[dg] (bus16)   at (5.5, 6.5)  [label=left:\scriptsize{$16$}] {};
		\node[bus] (bus16)   at (5.5, 6.5)   {};
		\node[bus] (bus15)   at (5.5, 7)  [label=left:\scriptsize{$15$}] {};
		\node[bus] (bus12)   at (6.5, 8)  [label=below:\scriptsize{$12$}] {};
		
		\node[substation] (bus11)   at (6.5, 8.5) {};  
		\node[bus] (bus11)   at (6.5, 8.5)   [label={[label distance=-.3mm]-175:\scriptsize{$11$}}] {};
		
		\node[bus] (bus10)   at (6.5, 9)  [label=left:\scriptsize{$10$}] {};
		\node[bus] (bus6)   at (2, 6.5)  [label=left:\scriptsize{$6$}] {};
		\node[bus] (bus5)   at (2, 7) [label=left:\scriptsize{$5$}] {};
		\node[bus] (bus4)   at (2, 7.5)  [label=left:\scriptsize{$4$}] {};
		\node[bus] (bus3)   at (2, 8)   [label={[label distance=-.3mm]-175:\scriptsize{$3$}}] {};
		\node[bus] (bus2)   at (2, 8.5)  [label=left:\scriptsize{$2$}] {};
		
		\node[substation] (bus1)   at (2, 9)  [label=right:\scriptsize{$1$}] {};
		\node[bus] (bus1)   at (2, 9)  {};
		
		\node[bus] (bus0)   at (2, 9.5)  [label=left:\scriptsize{$0$}] {};
		\node[bus] (bus18)   at (1, 9)  [label=below:\scriptsize{$18$}] {};
		\node[bus] (bus19)   at (0, 9)  [label=below:\scriptsize{$19$}] {};
		\node[bus] (bus22)   at (3, 8.5)  [label=above:\scriptsize{$22$}] {};
		\node[bus] (bus7)   at (4, 8.5)  [label=above:\scriptsize{$7$}] {};
		\node[bus] (bus23)   at (3, 8)  [label=right:\scriptsize{$23$}] {};
		\node[bus] (bus24)   at (3, 7.5)  [label=right:\scriptsize{$24$}] {};
		
		\node[substation] (bus25)   at (4, 7)  [label=right:\scriptsize{$25$}] {};
		\node[bus] (bus25)   at (4, 7)  {};
		
		\node[dg] (bus26)   at (4, 6.5)  [label=right:\scriptsize{$26$}] {};
		\node[bus] (bus26)   at (4, 6.5)  {};

		\node[bus] (bus27)   at (7.5, 9.5)  [label=right:\scriptsize{$27$}] {};
		\node[bus] (bus28)   at (7.5, 9)  [label=right:\scriptsize{$28$}] {};
		\node[bus] (bus29)   at (7.5, 8.5)  [label=right:\scriptsize{$29$}] {};
		\node[bus] (bus30)   at (7.5, 8)  [label=right:\scriptsize{$30$}] {};
		\node[bus] (bus9)   at (7.5, 7)  [label=right:\scriptsize{$9$}] {};
		
		\node[bus] (bus14)   at (1, 8)  [label={[label distance=-1mm]-30:\scriptsize{$14$}}] {};
		\node[bus] (bus8)   at (1, 7.5)  [label=below:\scriptsize{$8$}] {};
		\node[bus] (bus13)   at (0, 8)  [label=below:\scriptsize{$13$}] {};
		
		\node[bus] (bus21)   at (5.5, 8.5)  [label=left:\scriptsize{$21$}] {};
		\node[bus] (bus20)   at (5.5, 9)  [label=left:\scriptsize{$20$}] {};
		\draw [-] (bus17) -- (bus32);
		\draw [-] (bus32) -- (bus31);
		\draw [-] (bus17) -- (bus16);
		\draw [-] (bus16) -- (bus15);
		\draw [-] (bus15) -- (bus9);
		\draw [-] (bus9) -- (bus30);
		\draw [-] (bus30) -- (bus29);
		\draw [-] (bus29) -- (bus28);
		\draw [-] (bus28) -- (bus27);
		\draw [-] (bus11) -- (bus29);
		\draw [-] (bus11) -- (bus12);
		\draw [-] (bus11) -- (bus21);
		\draw [-] (bus11) -- (bus10);
		\draw [-] (bus21) -- (bus20);
		\draw [-] (bus25) -- (bus26);
		\draw [-] (bus25) -- (bus5);
		\draw [-] (bus5) -- (bus6);
		\draw [-] (bus5) -- (bus4);
		\draw [-] (bus4) -- (bus3);
		\draw [-] (bus3) -- (bus14);
		\draw [-] (bus14) -- (bus8);
		\draw [-] (bus14) -- (bus13);
		\draw [-] (bus1) -- (bus0);
		\draw [-] (bus1) -- (bus2);
		\draw [-] (bus1) -- (bus18);
		\draw [-] (bus18) -- (bus19);
		\draw [-] (bus2) -- (bus22);
		\draw [-] (bus22) -- (bus7);
		\draw [-] (bus22) -- (bus23);
		\draw [-] (bus23) -- (bus24);
		\draw [dashed,color=blue,line width=0.3mm] (4.8,9.8) -- (8.7,9.8);
		\draw [dashed,color=blue,line width=0.3mm] (8.7,9.8) -- (8.7,6.30);
		\draw [dashed,color=blue,line width=0.3mm] (8.7,6.30) -- (4.8,6.30);
		\draw [dashed,color=blue,line width=0.3mm] (4.8,6.30) -- (4.8,9.8);
		\draw [dashed,color=olive,line width=0.3mm] (4.5,6.30) -- (-0.5,6.30);
		\draw [dashed,color=olive,line width=0.3mm] (4.5,6.30) -- (4.5,7.20);
		\draw [dashed,color=olive,line width=0.3mm] (4.5,7.20) -- (2.3,7.20);
		\draw [dashed,color=olive,line width=0.3mm] (2.3,7.20) -- (2.3,8.20);
		\draw [dashed,color=olive,line width=0.3mm] (2.3,8.20) -- (-0.5,8.20);
		\draw [dashed,color=olive,line width=0.3mm] (-0.5,8.20) -- (-0.5,6.30);
		\draw [dashed,color=red,line width=0.3mm] (-0.5,8.35) -- (-0.5,9.80);
		\draw [dashed,color=red,line width=0.3mm] (-0.5,9.80) -- (4.5,9.80);
		\draw [dashed,color=red,line width=0.3mm] (4.5,9.80) -- (4.5,7.30);
		\draw [dashed,color=red,line width=0.3mm] (4.5,7.30) -- (2.5,7.30);
		\draw [dashed,color=red,line width=0.3mm] (2.5,7.30) -- (2.5,8.35);
		\draw [dashed,color=red,line width=0.3mm] (2.5,8.35) -- (-0.5,8.35);
		
		\node[legend1] (legend11)   at (9, 7.5)  [label=right:\scriptsize{$\mbox{DG position}$}] {};
		\node[legend2] (legend22)   at (9, 7)  [label=right:\scriptsize{$\mbox{Substation position}$}] {};
		\end{tikzpicture}
	}
	\end{center}
	\caption{Optimal configuration for \sred{the} 33-node distribution system} \label{33nodes}
\end{figure}

\begin{figure}[!h]
	\begin{center}
		\usepgflibrary{arrows}
		\usetikzlibrary{arrows}
		\usepgflibrary{decorations.pathmorphing} 
		\usetikzlibrary{decorations.pathmorphing} 
		\resizebox {6cm}{6cm}  {
		\begin{tikzpicture}
		[bus/.style={circle,inner sep=0pt,minimum size=4pt,fill=black},
		point/.style={circle,inner sep=0pt,minimum size=.1pt,fill=black},
		substation/.style={rectangle,draw=orange!100,thick,
			inner sep=0pt,minimum width=2.5mm,minimum height=2.5mm},
		dg/.style={regular polygon, regular polygon sides=3,inner sep=2pt ,draw=cyan!100, thick},legend1/.style={regular polygon, regular polygon sides=3,inner sep=1.5pt ,draw=cyan!100,fill=white, thick},legend2/.style={rectangle,draw=orange!100,fill=white!100,thick,
			inner sep=0pt,minimum width=2.5mm,minimum height=2.5mm}]
		\node[bus] (bus34)   at (1, 9)  [label=below:\scriptsize{$34$}] {};
		\node[bus] (bus33)   at (1.5, 9)  [label=below: \scriptsize{$33$}] {};
		\node[bus] (bus32)   at (2, 9)  [label=below:\scriptsize{$32$}] {};
		\node[bus] (bus31)   at (2.5, 9)  [label={[label distance=-.3mm]-45:\scriptsize{$31$}}] {};
		\node[bus] (bus30)   at (3, 9)  [label=above:\scriptsize{$30$}] {};
		\node[bus] (bus29)   at (3.5, 9)  [label=below:\scriptsize{$29$}] {};
		\node[bus] (bus28)   at (4, 9) [label=below:\scriptsize{$28$}] {};
		\node[bus] (bus27)   at (4.5, 9)   [label=below:\scriptsize{$27$}] {};
		\node[bus] (bus2)   at (5, 9)  [label=right:\scriptsize{$2$}] {};
		\node[substation] (bus1)   at (5, 9.5)  {};
		\node[bus] (bus1)   at (5, 9.5)  [label=right:\scriptsize{$1$}] {};
		\node[bus] (bus0)   at (5, 10)  [label=right:\scriptsize{$0$}] {};
		\node[dg] (bus3)   at (5, 8.5)  {};
		\node[bus] (bus3)   at (5, 8.5)  [label=left:\scriptsize{$3$}] {};
		\node[bus] (bus46)   at (6, 8.5)   [label={[label distance=-.3mm]-175:\scriptsize{$46$}}] {};
		\node[bus] (bus47)   at (6, 8)  [label=left:\scriptsize{$47$}] {};
		
		\node[bus] (bus59)   at (6, 7.5)  [label=left:\scriptsize{$59$}] {};
		\node[bus] (bus60)   at (6, 7)  [label=left:\scriptsize{$60$}] {};
		\node[bus] (bus4)   at (2.5, 8.5)  [label=right:\scriptsize{$4$}] {};
		\node[bus] (bus5)   at (2.5, 8)  [label=right:\scriptsize{$5$}] {};
		\node[bus] (bus6)   at (2.5, 7.5)  [label={[label distance=-.5mm]5:\scriptsize{$6$}}] {};
		\node[bus] (bus7)   at (2.5, 7)  [label=right:\scriptsize{$7$}] {};
		\node[bus] (bus8)   at (2.5, 6.5)  [label=right:\scriptsize{$8$}] {};
		\node[bus] (bus9)   at (2.5, 6)  [label=right:\scriptsize{$9$}] {};
		\node[bus] (bus51)   at (1.5, 6.5)  [label=left:\scriptsize{$51$}] {};
		\node[bus] (bus52)   at (1.5, 6)  [label=left:\scriptsize{$52$}] {};
		\node[bus] (bus63)   at (1.5, 5.5)  [label=left:\scriptsize{$63$}] {};
		\node[dg] (bus50)   at (1.5, 7) {};
		\node[bus] (bus50)   at (1.5, 7)  [label=right:\scriptsize{$50$}] {};
		\node[bus] (bus49)   at (1.5, 7.5)  [label=right:\scriptsize{$49$}] {};
		\node[bus] (bus58)   at (1, 7.5)  [label=below:\scriptsize{$58$}] {};
		\node[bus] (bus48)   at (1.5, 8)  [label=left:\scriptsize{$48$}] {};
		\node[bus] (bus55)   at (3.5, 7.5)  [label=above:\scriptsize{$55$}] {};
		\node[bus] (bus54)   at (3.5, 7)  [label=right:\scriptsize{$54$}] {};
		\node[bus] (bus56)   at (4, 7.5)  [label=above:\scriptsize{$56$}] {};
		\node[bus] (bus57)   at (4.5, 7.5)  [label=above:\scriptsize{$57$}] {};
		\node[point] (p1)   at (2, 9.5)   {};

		\node[bus] (bus68)   at (3, 5.5)  [label=left:\scriptsize{$68$}] {};
		\node[bus] (bus16)   at (3, 5)  [label={[label distance=-.3mm]-175:\scriptsize{$16$}}] {};
		\node[bus] (bus15)   at (3, 4.5)  [label=left:\scriptsize{$15$}] {};
		\node[bus] (bus53)   at (3, 4)  [label=below:\scriptsize{$53$}] {};
		\node[bus] (bus17)   at (2.5, 5)  [label=left:\scriptsize{$17$}] {};
		\node[bus] (bus67)   at (3.5, 5.5)  [label=above:\scriptsize{$67$}] {};
		\node[bus] (bus66)   at (3.5, 5)  [label=right:\scriptsize{$66$}] {};
		\node[bus] (bus65)   at (3.5, 4.5)  [label=below:\scriptsize{$65$}] {};
		\node[bus] (bus11)   at (5, 5.5)  [label=right:\scriptsize{$11$}] {};
		\node[bus] (bus10)   at (5, 6)  [label=right:\scriptsize{$10$}] {};
		\node[bus] (bus42)   at (5, 6.5)  [label=right:\scriptsize{$42$}] {};
		\node[bus] (bus12)   at (5, 5)  [label={[label distance=-.3mm]-175:\scriptsize{$12$}}] {};
		\node[substation] (bus13)   at (5, 4.5)  {};
		\node[bus] (bus13)   at (5, 4.5)  [label=left:\scriptsize{$13$}] {};
		\node[bus] (bus24)   at (6, 4.5)  [label=right:\scriptsize{$24$}] {};
		\node[bus] (bus23)   at (6, 5)  [label=right:\scriptsize{$23$}] {};
		\node[bus] (bus22)   at (6, 5.5)  [label=right:\scriptsize{$22$}] {};
		\node[bus] (bus21)   at (6, 6)  [label=right:\scriptsize{$21$}] {};
		\node[dg] (bus25)   at (6, 4) {};
		\node[bus] (bus25)   at (6, 4)  [label=right:\scriptsize{$25$}] {};
		\node[bus] (bus26)   at (6, 3.5)  [label=right:\scriptsize{$26$}] {};
		\node[bus] (bus14)   at (5, 4)  [label=below:\scriptsize{$14$}] {};
		\node[bus] (bus45)   at (5.5, 4)  [label=below:\scriptsize{$45$}] {};
		\node[bus] (bus18)   at (5, 2.5)  [label=right:\scriptsize{$18$}] {};
		\node[bus] (bus19)   at (5, 2)  [label=right:\scriptsize{$19$}] {};
		\node[dg] (bus20)   at (5, 1.5)  {};
		\node[bus] (bus20)   at (5, 1.5)  [label=right:\scriptsize{$20$}] {};
		\node[point] (p2)   at (4, 5)   {};
		\node[point] (p3)   at (4, 1.5)  {};

		\node[bus] (bus43)   at (6, 9.5)  [label=above:\scriptsize{$43$}] {};
		\node[bus] (bus44)   at (7, 9.5)  [label=right:\scriptsize{$44$}] {};
		\node[bus] (bus35)   at (7, 8.5)  [label=right:\scriptsize{$35$}] {};
		\node[bus] (bus36)   at (7, 8)  [label=right:\scriptsize{$36$}] {};
		\node[bus] (bus37)   at (7, 7.5)  [label=right:\scriptsize{$37$}] {};
		\node[bus] (bus38)   at (7, 7)  [label=right:\scriptsize{$38$}] {};
		\node[substation] (bus39)   at (7, 6.5)  {};
		\node[bus] (bus39)   at (7, 6.5)  [label=right:\scriptsize{$39$}] {};
		\node[bus] (bus40)   at (7, 6)  [label=right:\scriptsize{$40$}] {};
		\node[bus] (bus41)   at (7, 5.5)  [label=right:\scriptsize{$41$}] {};

		\node[bus] (bus64)   at (1.5, 4.5)  [label=left:\scriptsize{$64$}] {};
		\node[substation] (bus61)   at (1.5, 4)  {};
		\node[bus] (bus61)   at (1.5, 4)  [label=left:\scriptsize{$61$}] {};
		\node[bus] (bus62)   at (1.5, 3.5)  [label=left:\scriptsize{$62$}] {};
		\draw [-] (bus0) -- (bus1);
		\draw [-] (bus1) -- (bus2);
		\draw [-] (bus2) -- (bus3);
		\draw [-] (bus3) -- (bus46);
		\draw [-] (bus46) -- (bus47);
		\draw [-] (bus47) -- (bus59);
		\draw [-] (bus59) -- (bus60);
		\draw [-] (bus2) -- (bus27);
		\draw [-] (bus27) -- (bus28);
		\draw [-] (bus28) -- (bus29);
		\draw [-] (bus30) -- (bus31);
		\draw [-] (bus31) -- (bus32);
		\draw [-] (bus32) -- (bus33);
		\draw [-] (bus33) -- (bus34);
		\draw [-] (bus1) -- (p1) -- (bus32);
		\draw [-] (bus31) -- (bus4);
		\draw [-] (bus4) -- (bus5);
		\draw [-] (bus5) -- (bus6);
		\draw [-] (bus6) -- (bus7);
		\draw [-] (bus7) -- (bus8);
		\draw [-] (bus8) -- (bus9);
		\draw [-] (bus6) -- (bus55);
		\draw [-] (bus55) -- (bus54);
		\draw [-] (bus55) -- (bus56);
		\draw [-] (bus56) -- (bus57);
		\draw [-] (bus8) -- (bus51);
		\draw [-] (bus51) -- (bus52);
		\draw [-] (bus52) -- (bus63);
		\draw [-] (bus51) -- (bus50);
		\draw [-] (bus50) -- (bus49);
		\draw [-] (bus49) -- (bus58);
		\draw [-] (bus49) -- (bus48);
		
		\draw [-] (bus43) -- (bus44);
		\draw [-] (bus44) -- (bus35);
		\draw [-] (bus35) -- (bus36);
		\draw [-] (bus36) -- (bus37);
		\draw [-] (bus37) -- (bus38);
		\draw [-] (bus38) -- (bus39);
		\draw [-] (bus39) -- (bus40);
		\draw [-] (bus40) -- (bus41);
		
		\draw [-] (bus68) -- (bus67);
		\draw [-] (bus67) -- (bus66);
		\draw [-] (bus66) -- (bus65);
		\draw [-] (bus67) -- (bus11);
		\draw [-] (bus11) -- (bus10);
		\draw [-] (bus10) -- (bus42);
		\draw [-] (bus11) -- (bus12);
		\draw [-] (bus12) -- (bus13);
		\draw [-] (bus13) -- (bus24);
		\draw [-] (bus24) -- (bus23);
		\draw [-] (bus23) -- (bus22);
		\draw [-] (bus22) -- (bus21);
		\draw [-] (bus24) -- (bus25);
		\draw [-] (bus25) -- (bus26);
		\draw [-] (bus13) -- (bus14);
		\draw [-] (bus14) -- (bus45);
		\draw [-] (bus15) -- (bus53);
		\draw [-] (bus15) -- (bus16);
		\draw [-] (bus16) -- (bus17);
		\draw [-] (bus18) -- (bus19);
		\draw [-] (bus19) -- (bus20);
		\draw [-] (bus16) -- (bus68);
		\draw [-] (bus12) -- (p2) -- (p3) -- (bus20);

		\draw [-] (bus64) -- (bus61);
		\draw [-] (bus61) -- (bus62);
		
		\node[point] (p10)   at (.75, 10.3)  {};
		\node[point] (p11)   at (5.4, 10.3)  {};
		\node[point] (p12)   at (5.4, 9)  {};
		\node[point] (p13)   at (6.4, 9)  {};
		\node[point] (p14)   at (6.4, 6.8)  {};
		\node[point] (p15)   at (3, 6.8)  {};
		\node[point] (p16)   at (3, 5.8)  {};
		\node[point] (p17)   at (1.7, 5.8)  {};
		\node[point] (p18)   at (1.7, 5.3)  {};
		\node[point] (p19)   at (.75, 5.3)  {};
		
		\draw [dashed,color=red,line width=0.3mm] (p10) -- (p11);	
		\draw [dashed,color=red,line width=0.4mm] (p11) -- (p12);
		\draw [dashed,color=red,line width=0.3mm] (p12) -- (p13);
		\draw [dashed,color=red,line width=0.4mm] (p13) -- (p14);
		\draw [dashed,color=red,line width=0.3mm] (p14) -- (p15);
		\draw [dashed,color=red,line width=0.4mm] (p15) -- (p16);
		\draw [dashed,color=red,line width=0.3mm] (p16) -- (p17);
		\draw [dashed,color=red,line width=0.4mm] (p17) -- (p18);
		\draw [dashed,color=red,line width=0.3mm] (p18) -- (p19);
		\draw [dashed,color=red,line width=0.3mm] (p19) -- (p10);
		
		\node[point] (p24)   at (2, 1)  {};
		\node[point] (p25)   at (2, 5.7)  {};
		\node[point] (p26)   at (3.2, 5.7)  {};
		\node[point] (p27)   at (3.2, 6.7)  {};
		\node[point] (p28)   at (6.5, 6.7)  {};
		\node[point] (p29)   at (6.5, 3.2)  {};
		\node[point] (p30)   at (5.5, 3.2)  {};
		\node[point] (p31)   at (5.5, 1)  {};
		\node[point] (p32)   at (2, 1)  {};

		\draw [dashed,color=olive,line width=0.4mm] (p24) -- (p25);
		\draw [dashed,color=olive,line width=0.3mm] (p25) -- (p26);
		\draw [dashed,color=olive,line width=0.4mm] (p26) -- (p27);
		\draw [dashed,color=olive,line width=0.3mm] (p27) -- (p28);
		\draw [dashed,color=olive,line width=0.4mm] (p28) -- (p29);
		\draw [dashed,color=olive,line width=0.3mm] (p29) -- (p30);
		\draw [dashed,color=olive,line width=0.4mm] (p30) -- (p31);
		\draw [dashed,color=olive,line width=0.4mm] (p31) -- (p32);
		
		\node[point] (p33)   at (5.7, 10.3)  {};
		\node[point] (p34)   at (5.7, 9.25)  {};
		\node[point] (p35)   at (6.7, 9.25)  {};
		\node[point] (p36)   at (6.7, 5)  {};
		\node[point] (p37)   at (7.6, 5)  {};
		\node[point] (p38)   at (7.6, 10.3)  {};
		
		\draw [dashed,color=blue,line width=0.4mm] (p33) -- (p34);
		\draw [dashed,color=blue,line width=0.4mm] (p34) -- (p35);
		\draw [dashed,color=blue,line width=0.4mm] (p35) -- (p36);
		\draw [dashed,color=blue,line width=0.4mm] (p36) -- (p37);
		\draw [dashed,color=blue,line width=0.4mm] (p37) -- (p38);
		\draw [dashed,color=blue,line width=0.4mm] (p38) -- (p33);

		\node[point] (p20)   at (.75, 4.9)  {};
		\node[point] (p21)   at (1.75, 4.9)  {};
		\node[point] (p22)   at (1.75, 3)  {};
		\node[point] (p23)   at (.75, 3)  {};
		
		\draw [dashed,color=purple,line width=0.3mm] (p20) -- (p21);
		\draw [dashed,color=purple,line width=0.3mm] (p21) -- (p22);
		\draw [dashed,color=purple,line width=0.3mm] (p22) -- (p23);
		\draw [dashed,color=purple,line width=0.3mm] (p23) -- (p20);

		\node[legend1] (legend11)   at (6, 2)  [label=right:\scriptsize{$\mbox{DG position}$}] {};
		\node[legend2] (legend22)   at (6, 1.5)  [label=right:\scriptsize{$\mbox{Substation position}$}] {};
		\end{tikzpicture}
	}
	\end{center}
	\caption{Optimal configuration for \sred{the} 69-node distribution system} \label{69nodes}
\end{figure}
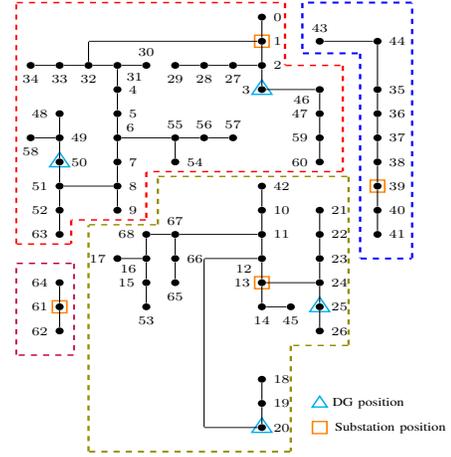
\subsection{Optimal distribution network configuration}
\sred{We report optimal configurations for} the 33-node and \sred{the} 69-node distribution systems in Fig. \ref{33nodes} and Fig. \ref{69nodes}, respectively.
\begin{table}[h!]
	\caption{\sred{Comparison of load shedding}} \label{comparing}
	\centering
	\setlength{\extrarowheight}{1.4pt}
	\begin{tabularx}{0.5\textwidth}{X X X X X X X X X}
		\hline
		$ $ & \multicolumn{4}{c}{\centering {DR model}} & \multicolumn{4}{c}{\centering {Robust model}}\\
		\cmidrule(r{1em}){2-5}
		\cmidrule(r{0em}){6-9}
		 $ $ & WCD & WCS & Sim & Time & WCD & WCS & Sim & Time \\
		Nodes & (KW) & (KW) & (KW) & (s) & (KW) & (KW) & (KW) & (s) \\
		\hline
		 33 & 1655 & 2535 & 1451 & 91 & 1921 & 2352 & 1648 & 179 \\
		 69 & 4297 & 5119 & 3594 & 173 & 4570 & 4998 & 4014 & 372\\
		\hline
	\end{tabularx}
\end{table}
\sred{We also compare our DR model with the RO model.}
For comparison purposes, we fix configuration decisions obtained by each model and then simulate the load shedding using randomly generated contingencies. Table \ref{comparing} reports the expected load shedding under the worst-case contingency distribution (WCD), the load shedding under the worst-case contingency scenario (WCS), the out-of-sample average load shedding under a randomly simulated contingency distribution within $\mathscr{D}$ (Sim), and the computational time of both models. The results verify that our DR approach yields lower load shedding both under worst-case distribution and in the out-of-sample simulations. In particular, our approach leads to 11\% and 10\% reduction in average load shedding in the out-of-sample simulation and 13\% and 5\% reduction under the worst-case distribution for the 33-node and 69-node distribution systems, respectively. For the worst-case contingency scenario, RO model triggers less load shedding, which was expected, because RO optimizes the system configuration with respect to the worst-case contingency scenario. In addition, the CPU seconds taken to solve the test instances demonstrate the efficacy of the proposed solution approach. To further verify the efficacy, we replicate the experiments on 10 randomly generated instances. For the 33-node system, the average and maximum number of iterations the CCG algorithm takes to converge are 7.6 and 12, respectively; and for the 69-node system, the average and maximum number of iterations are 8.4 and 14, respectively.

\begin{figure}[h!]
	\centering
	\includegraphics[width=2.3in]{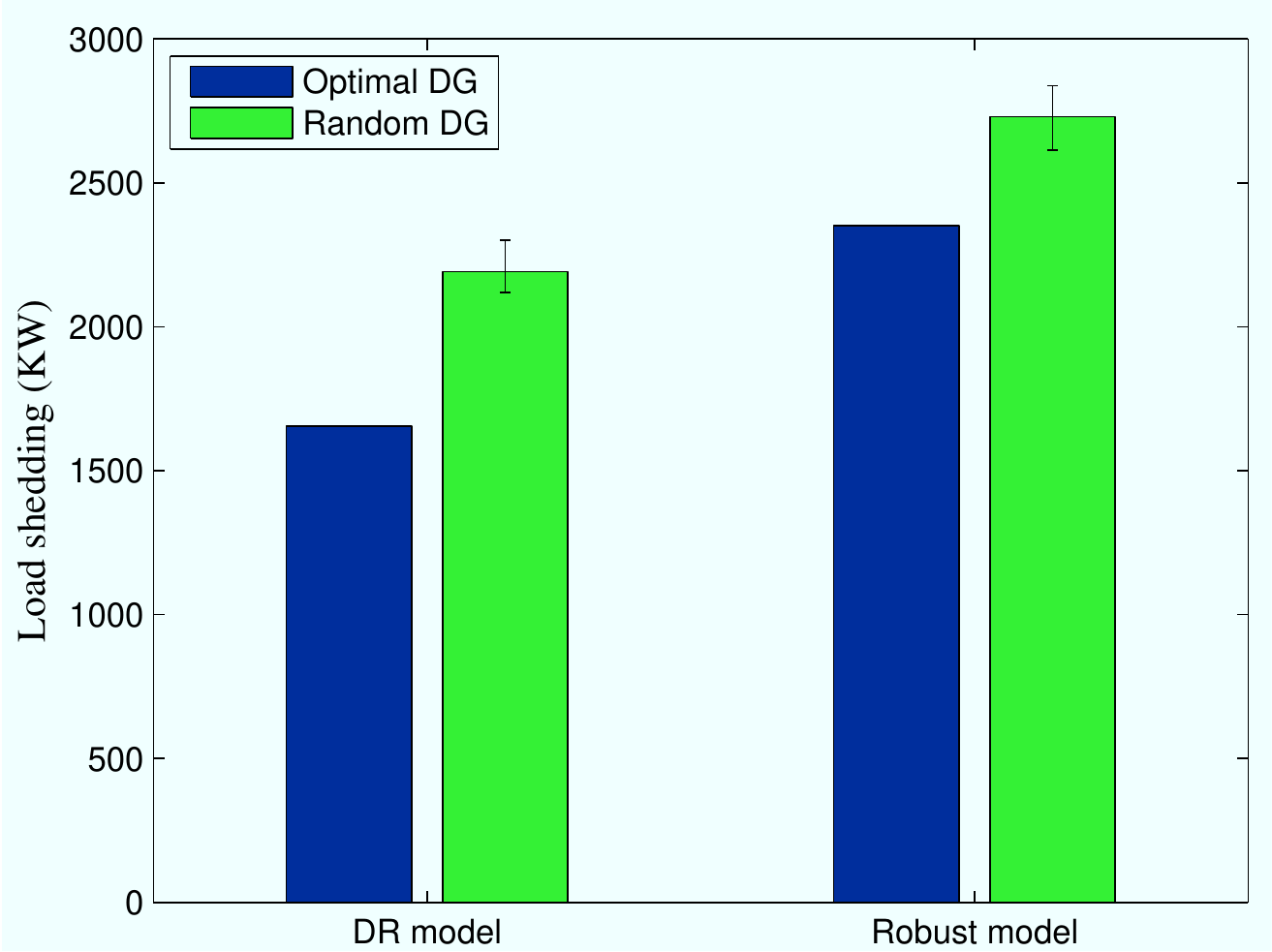}
	\caption{Comparisons of optimal and random DG allocation in the 33-node distribution system} \label{figH33}
\end{figure}

\begin{figure}[h!]
	\centering
	\includegraphics[width=2.3in]{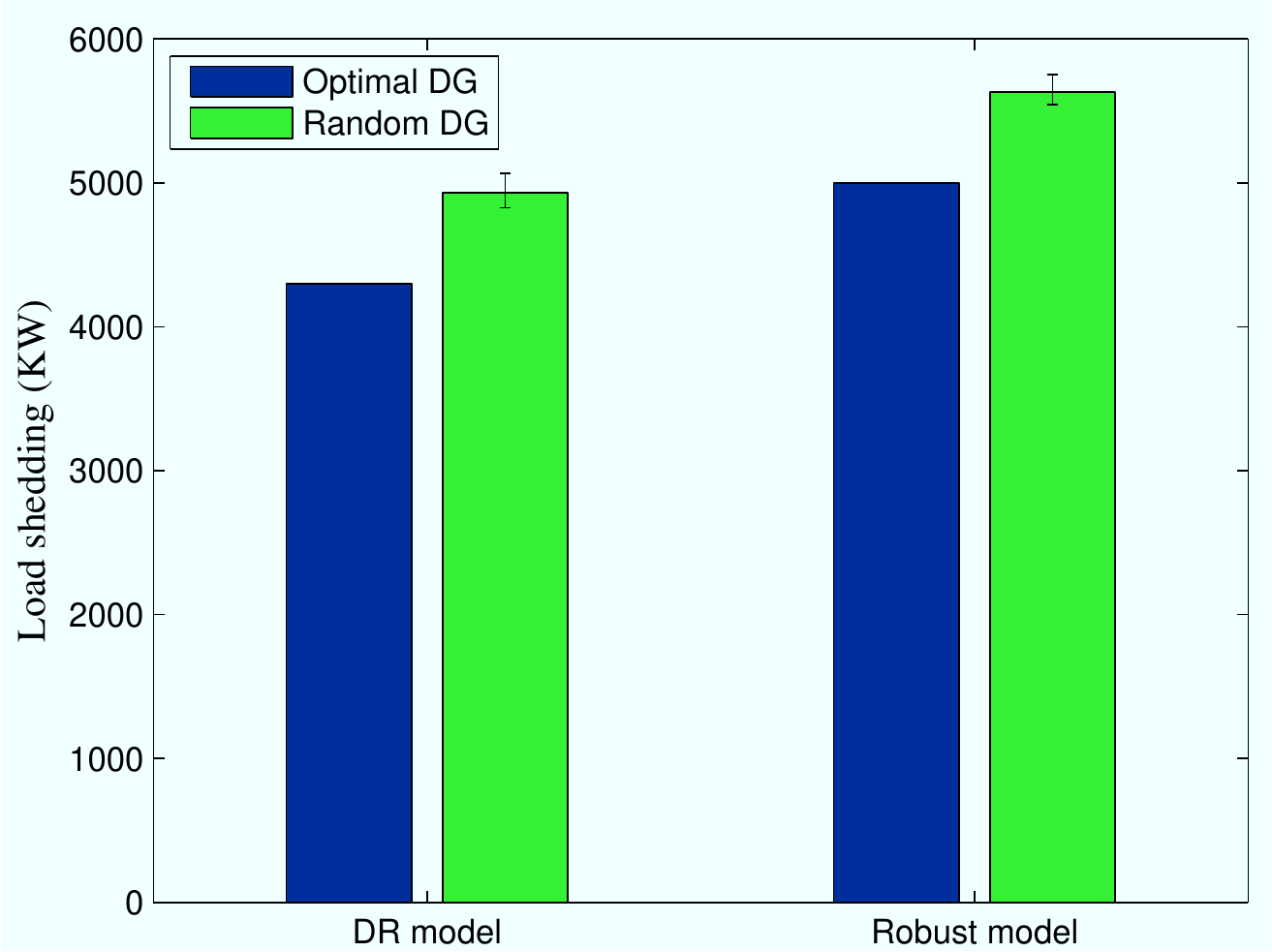}
	\caption{Comparisons of optimal and random DG allocation in the 69-node distribution system} \label{figH69}
\end{figure}
\subsection{On the value of optimal DG allocation}
We conduct a set of experiments to evaluate the value of optimally allocating DG units in the distribution system. In Fig. \ref{figH33} and Fig. \ref{figH69} we compare the level of load shedding when DG units are optimally located with the case when DG units are randomly deployed. For ``optimal DG", we solve the DR and RO models. For ``random DG", we first randomly place DGs and then solve both models to configure the distribution system. We perform the experiments for 5 times and report the average values to mitigate the randomness.
From Figs. \ref{figH33} and \ref{figH69}, we observe that locating DGs properly can significantly decrease the load shedding. This is because when the distribution system is affected by contingencies, the loads in islanded zones can be effectively picked up by the existing DG resources. As a result, better DG allocation significantly enhances the system resiliency.

\begin{figure}[h!]
	\centering
	\includegraphics[width=6cm,height=4cm]{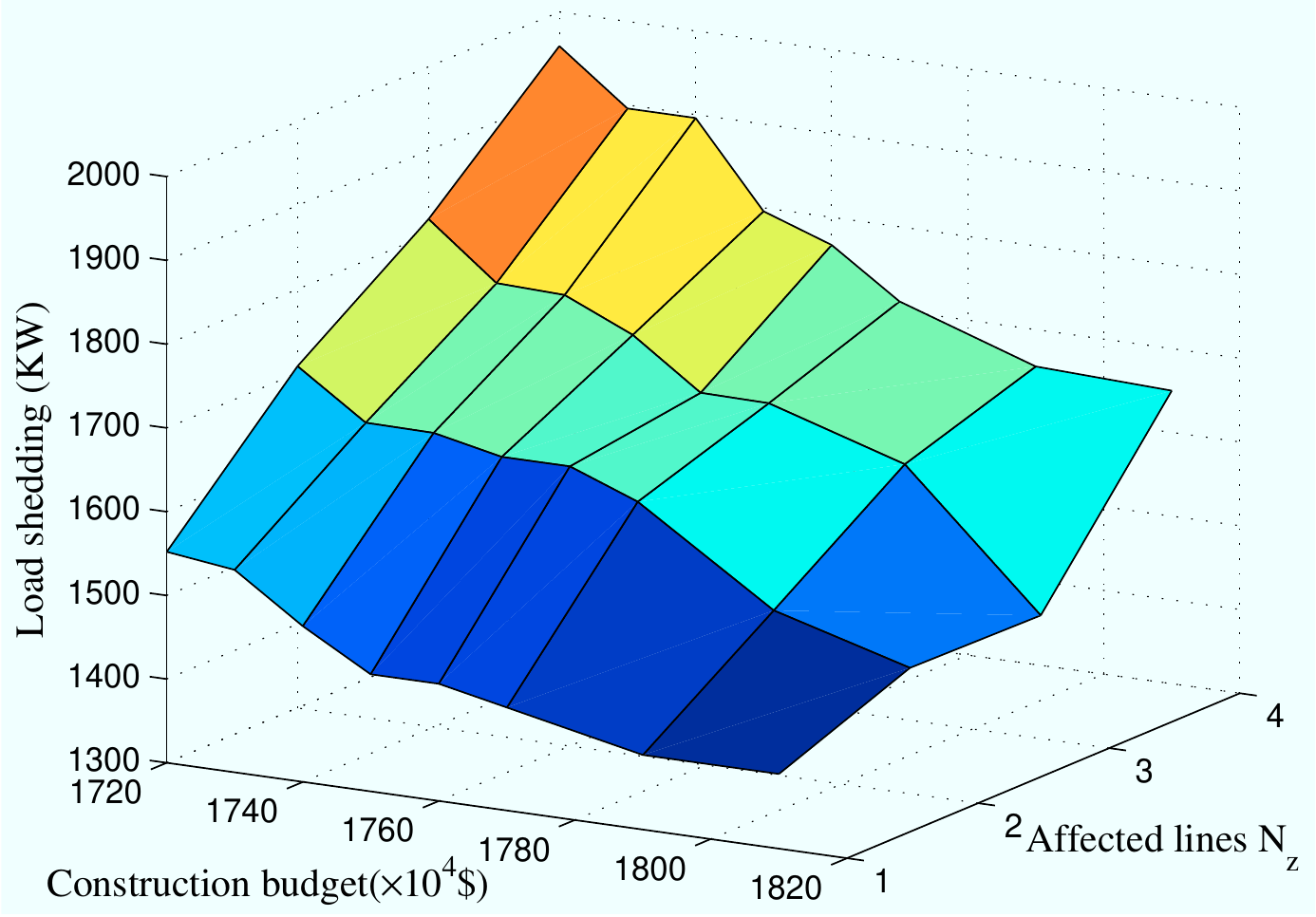}
	\caption{Average load shedding under various line construction budget and affected lines for the 33-node distribution system} \label{fig33D}
\end{figure}

\begin{figure}[h!]
	\centering
	\includegraphics[width=6cm,height=4cm]{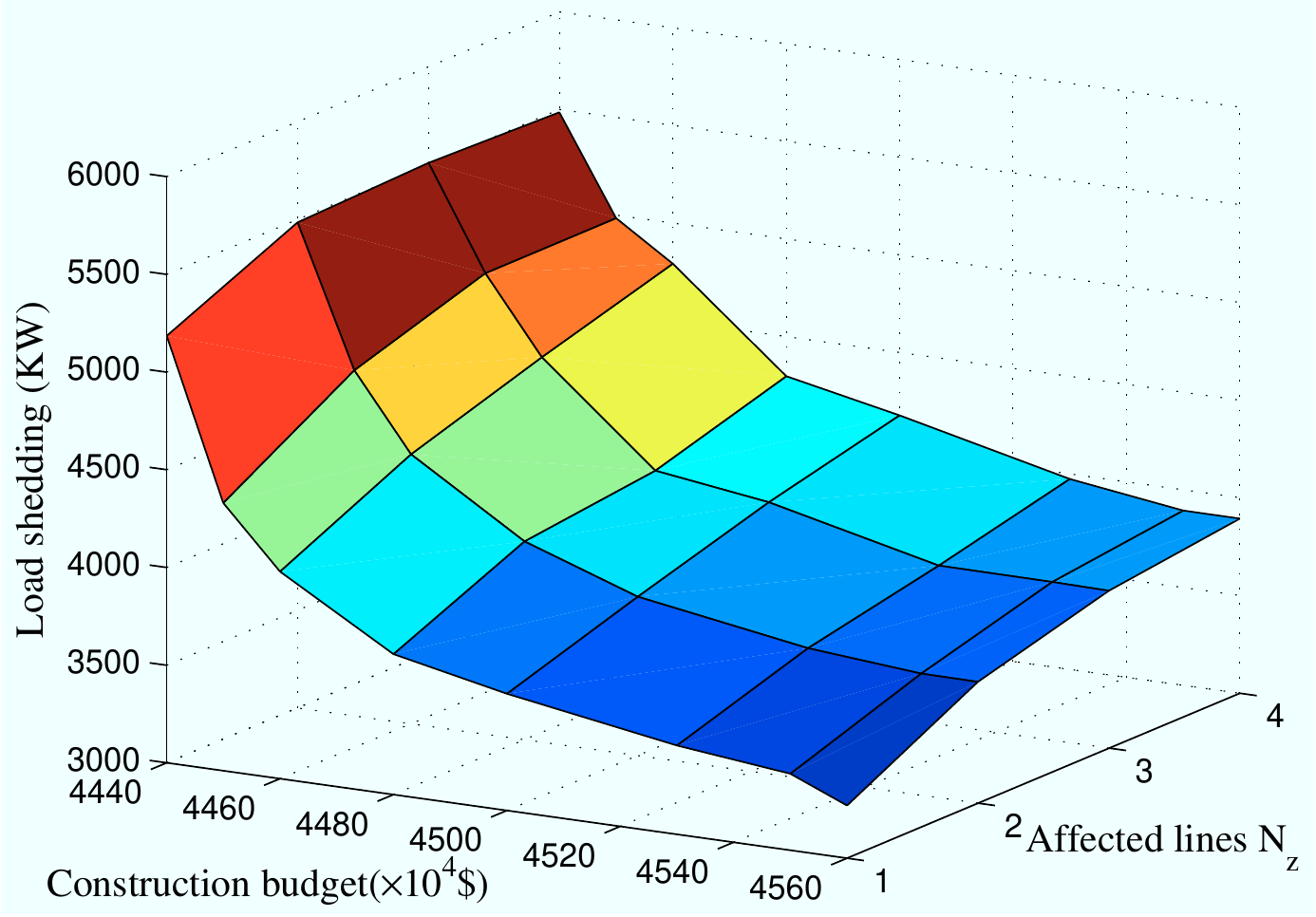}
	\caption{Average load shedding under various line construction budget and affected lines for the 69-node distribution system} \label{fig69D}
\end{figure}

\subsection{Impact of construction and contingency budgets}
{\color{black}In Figs. \ref{fig33D} and \ref{fig69D}, we depict the amounts of expected load shedding under various line construction budgets (i.e., $B_y$) and contingency budgets (i.e., $N_z$). From these two figures, we observe that load shedding reduces as $B_y$ increases and as $N_z$ decreases, i.e., as we allow the contingency to affect less power lines in the DR model. This is intuitive. In addition, we observe that load shedding is sensitive to the construction budget. For example, by increasing the budget from \$4440$\times 10^4$ to \$4480$\times 10^4$ when $N_z = 4$ in Fig. 7, the load shedding decreases from 5485KW to 4297KW, which means that a 0.9\% budgetary rise translates into a 21.6\% load shedding reduction. Furthermore, we observe that the impact of construction budget is marginally diminishing. For example, increasing the budget from \$4500$\times 10^4$ to \$4560$\times 10^4$ (i.e., by 1.3\%) results in a 6.7\% load shedding reduction. This observation highlights the necessity of implementing a cost-effective distribution configuration planning. }
\subsection{Worst-case contingency distribution}
{\color{black}The worst-case contingency distribution for the 69-node distribution system is reported in Table \ref{table2}.
{\color{black} We select a subset of representative scenarios to display and omit other scenarios with smaller probability values.}
From this table, we observe that the contingency probabilities for different power lines are highly heterogeneous. This provides the system operator a guideline on the system vulnerability and a meaningful contingency probability distribution that can be used in other vulnerability analyses.}
\begin{table}[h!]
	\caption{Worst-case contingency distribution for the 69-node system} \label{table2}
	\centering
	\begin{tabular}{c c c}
		\hline Scenario & Affected lines & Probability\\
		\hline
		1 & 6-15,13-14,34-35,39-40 & 0.0031\\
		2 & 5-26,6-15,12-13,38-39 & 0.0025 \\
		3 & 1-2,5-26,38-39,39-40 & 0.0021\\
		4 & 12-13,13-14,38-39,61-62 & 0.0002\\
		\hline
	\end{tabular}
\end{table}

\appendix{}
\begin{IEEEproof} [\hspace{-.3in} Proof of Proposition \ref{prop:ref}]
We rewrite $\max_{\mathbb{P} \in \mathbb{D}} E_{\mathbb{P}}[Q(\mathbf{g},{\mathbf{z}})]$ as:
\begin{subequations}
\begin{eqnarray}
&&\hspace{-0.5in} \max_{\mathbb{P} \in \mathbb{D}} E_{\mathbb{P}}[Q(\mathbf{g},{\mathbf{z}})]= \max_{\mathbb{P}} \int_{\mathfrak{D{(\mathbf{g})}}}{Q(\mathbf{g},\mathbf{z})d{\mathbb{P}}}, \label{dual1} \\
&\hspace{-0.5in} \mbox{s.t.} & \hspace{-0.1in}\int_{\mathfrak{D{(\mathbf{g})}}}{d{\mathbb{P}}}=1, \label{dual2}  \\
&&\hspace{-0.5in} \int_{\mathfrak{D{(\mathbf{g})}}}{(1-z_{mn,t}) d{\mathbb{P}}} \leq \mu_{mn,t}^{max}, \forall{(m,n) \in \mathcal{E}}, \  \forall{t \in \mathcal{T}}. \label{dual3}
\end{eqnarray}
\end{subequations}
The feasible region of \sred{the} problem \eqref{dual1}--\eqref{dual3} has an interior point. In other words, there exists \sred{a} $\hat{\mathbb{P}}$ \sred{that satisfies constraint \eqref{dual2} at equality and constraint \eqref{dual3} strictly. }
\sred{For example, we can set $\hat{\mathbb{P}}$ to be the probability distribution solely supported on the scenario that}
no contingency arises in the system, i.e., $z_{mn,t}=1, \forall{(m,n) \in \mathcal{E}, t \in \mathcal{T}}$. Thus, \sred{the} Slater's condition holds
\sred{between the problem \eqref{dual1}--\eqref{dual3} and the following dual formulation:}
	\begin{eqnarray}
	&&\hspace{-0.5in} \min_{\mathbf{\boldsymbol\beta} \geq 0,\gamma} \gamma + \sum_{t \in \mathcal{T}} \sum_{(m,n) \in \mathcal{E}} \mu_{mn,t}^{max}  \beta_{mn,t}, \label{dual4} \\
	&&\hspace{-0.5in} \mbox{s.t.} \nonumber\\
	&&\hspace{-0.5in} \gamma + \sum_{t \in \mathcal{T}} \sum_{(m,n) \in \mathcal{E}} (1-z_{mn,t}) \beta_{mn,t} \geq  Q(\mathbf{g},\mathbf{z}), \ \forall{\mathbf{z} \in \mathfrak{D{(\mathbf{g})}}}. \label{dual5}
	\end{eqnarray}
	where $\gamma$ and $\mathbf{\boldsymbol\beta}$ are dual variables \sred{associated with} constraints \eqref{dual2} and \eqref{dual3}, respectively. In \sred{the dual} formulation, we observe that the optimal $\gamma$ should satisfy
	\begin{eqnarray}
	&&\hspace{-0.5in} \gamma = \max_{\mathbf{z} \in \mathfrak{D{(\mathbf{g})}}} \bigg\{Q(\mathbf{g},\mathbf{z})- \sum_{t \in \mathcal{T}} \sum_{(m,n) \in \mathcal{E}} (1-z_{mn,t}) \beta_{mn,t}\bigg\}. \label{gamma}
	\end{eqnarray}
	Substituting $\gamma$ from \eqref{gamma} to \sred{the} objective function \eqref{dual4} \sred{completes} the proof.
\end{IEEEproof}

\vspace{.3in}

\begin{IEEEproof} [\hspace{-.3in} Proof of Proposition \ref{prop:wc-distribution}]
With variables $g$ and $u^j$ fixed at $\hat{g}^R$ and $\hat{u}^j$, respectively, we take the dual of formulation \eqref{master1}--\eqref{master3} to obtain:
\begin{subequations}
\begin{align}
\hspace{-0.5cm} \max_{\psi \geq 0} \ & \ \sum_{j=1}^R \psi_j \Bigl( \sum_{t \in \mathcal{T}} \sum_{n \in \mathcal{N}} s^j_{nt} \Bigr) \label{dual-obj} \\
\hspace{-0.5cm} \mbox{s.t.} \ & \ \sum_{j=1}^R \psi_j (1 - z^j_{mn, t}) \leq \mu^{max}_{mn, t}, \nonumber \\
& \ \forall t \in \mathcal{T}, \ \ \forall (m, n) \in \mathcal{E}, \label{dual-con-1} \\
\hspace{-0.5cm} & \ \sum_{j=1}^R \psi_j = 1. \label{dual-con-2}
\end{align}
\end{subequations}
By constraints \eqref{dual-con-1}--\eqref{dual-con-2}, $\{\psi_j\}_{j = 1, \ldots, R}$ characterize a probability distribution supported on scenarios $\{\mathbf{z}^j\}_{j = 1, \ldots, R}$ such that $\mathbb{P}\{\mathbf{z} = \mathbf{z}^j\} = \psi_j$, $\forall j = 1, \ldots, R$. As the CCG framework terminates at the $R^{\mbox{\tiny th}}$ iteration and by the strong duality of linear programming, formulation \eqref{dual-obj}--\eqref{dual-con-2} is equivalent to the worst-case expectation formulation \eqref{worstcaseproblem}, i.e., these two formulations yield the same optimal value. It follows that $\{\psi_j\}_{j = 1, \ldots, R}$ characterize the worst-case contingency probability distribution.
\end{IEEEproof}

\bibliographystyle{IEEEtran}

\begin{thebibliography}{10}
\providecommand{\url}[1]{#1}
\csname url@samestyle\endcsname
\providecommand{\newblock}{\relax}
\providecommand{\bibinfo}[2]{#2}
\providecommand{\BIBentrySTDinterwordspacing}{\spaceskip=0pt\relax}
\providecommand{\BIBentryALTinterwordstretchfactor}{4}
\providecommand{\BIBentryALTinterwordspacing}{\spaceskip=\fontdimen2\font plus
\BIBentryALTinterwordstretchfactor\fontdimen3\font minus
  \fontdimen4\font\relax}
\providecommand{\BIBforeignlanguage}[2]{{%
\expandafter\ifx\csname l@#1\endcsname\relax
\typeout{** WARNING: IEEEtran.bst: No hyphenation pattern has been}%
\typeout{** loaded for the language `#1'. Using the pattern for}%
\typeout{** the default language instead.}%
\else
\language=\csname l@#1\endcsname
\fi
#2}}
\providecommand{\BIBdecl}{\relax}
\BIBdecl

\bibitem{salman2015evaluating}
A.~M. Salman, Y.~Li, and M.~G. Stewart, ``Evaluating system reliability and
  targeted hardening strategies of power distribution systems subjected to
  hurricanes,'' \emph{Reliability Engineering \& System Safety}, vol. 144, pp.
  319--333, 2015.

\bibitem{che2014only}
L.~Che, M.~Khodayar, and M.~Shahidehpour, ``Only connect: Microgrids for
  distribution system restoration,'' \emph{IEEE Power and Energy Magazine},
  vol.~12, no.~1, pp. 70--81, 2014.

\bibitem{ofeconomic}
E.~O. of~the President, ``Economic benefits of increasing electric grid
  resilience to weather outages-august 2013.''

\bibitem{ton2015more}
D.~T. Ton and W.~P. Wang, ``A more resilient grid: The {U.S. Department of
  Energy} joins with stakeholders in an {R\&D} plan,'' \emph{IEEE Power and
  Energy Magazine}, vol.~13, no.~3, pp. 26--34, 2015.

\bibitem{report}
``{President's Council of Economic Advisers and the U.S. Department of Energy,
  Economic Benefits of Increasing Electric Grid Resilience to Weather
  Outages},'' Tech. Rep., 2013.

\bibitem{moreira2011large}
J.~Moreira, E.~Miguez, C.~Vilacha, and A.~F. Otero, ``Large-scale network
  layout optimization for radial distribution networks by parallel computing,''
  \emph{IEEE transactions on power delivery}, vol.~26, no.~3, pp. 1946--1951,
  2011.

\bibitem{kumar2014design}
D.~Kumar and S.~Samantaray, ``Design of an advanced electric power distribution
  systems using seeker optimization algorithm,'' \emph{International Journal of
  Electrical Power \& Energy Systems}, vol.~63, pp. 196--217, 2014.

\bibitem{wang2015self}
Z.~Wang and J.~Wang, ``Self-healing resilient distribution systems based on
  sectionalization into microgrids,'' \emph{IEEE Transactions on Power
  Systems}, vol.~30, no.~6, pp. 3139--3149, 2015.

\bibitem{arefifar2013comprehensive}
S.~A. Arefifar, Y.~A.-R.~I. Mohamed, and T.~H. EL-Fouly, ``Comprehensive
  operational planning framework for self-healing control actions in smart
  distribution grids,'' \emph{IEEE Transactions on Power Systems}, vol.~28,
  no.~4, pp. 4192--4200, 2013.

\bibitem{yuan2016robust}
W.~Yuan, J.~Wang, F.~Qiu, C.~Chen, C.~Kang, and B.~Zeng, ``Robust
  optimization-based resilient distribution network planning against natural
  disasters,'' \emph{IEEE Transactions on Smart Grid}, vol.~7, no.~6, pp.
  2817--2826, 2016.

\bibitem{ma2016resilience}
S.~Ma, B.~Chen, and Z.~Wang, ``Resilience enhancement strategy for distribution
  systems under extreme weather events,'' \emph{IEEE Transactions on Smart
  Grid}, vol.~9, no.~2, pp. 1442--1451, 2016.

\bibitem{samui2012direct}
A.~Samui, S.~Singh, T.~Ghose, and S.~Samantaray, ``A direct approach to optimal
  feeder routing for radial distribution system,'' \emph{IEEE Transactions on
  Power Delivery}, vol.~27, no.~1, pp. 253--260, 2012.

\bibitem{navarro2009large}
A.~Navarro and H.~Rudnick, ``Large-scale distribution planning--part {I}:
  Simultaneous network and transformer optimization,'' \emph{IEEE Transactions
  on Power Systems}, vol.~24, no.~2, pp. 744--751, 2009.

\bibitem{mahmoud2016optimal}
K.~Mahmoud, N.~Yorino, and A.~Ahmed, ``Optimal distributed generation
  allocation in distribution systems for loss minimization,'' \emph{IEEE
  Transactions on Power Systems}, vol.~31, no.~2, pp. 960--969, 2016.

\bibitem{pereira2016optimal}
B.~R. Pereira, G.~R.~M. da~Costa, J.~Contreras, and J.~R.~S. Mantovani,
  ``Optimal distributed generation and reactive power allocation in electrical
  distribution systems,'' \emph{IEEE Transactions on Sustainable Energy},
  vol.~7, no.~3, pp. 975--984, 2016.

\bibitem{atwa2010optimal}
Y.~Atwa, E.~El-Saadany, M.~Salama, and R.~Seethapathy, ``Optimal renewable
  resources mix for distribution system energy loss minimization,'' \emph{IEEE
  Transactions on Power Systems}, vol.~25, no.~1, pp. 360--370, 2010.

\bibitem{liu2011optimal}
Z.~Liu, F.~Wen, and G.~Ledwich, ``Optimal siting and sizing of distributed
  generators in distribution systems considering uncertainties,'' \emph{IEEE
  Transactions on power delivery}, vol.~26, no.~4, pp. 2541--2551, 2011.

\bibitem{miguez2002improved}
E.~M{\'\i}guez, J.~Cidr{\'a}s, E.~D{\'\i}az-Dorado, and J.~L.
  Garc{\'\i}a-Dornelas, ``An improved branch-exchange algorithm for large-scale
  distribution network planning,'' \emph{IEEE Transactions on Power Systems},
  vol.~17, no.~4, pp. 931--936, 2002.

\bibitem{abdelaziz2010distribution}
A.~Y. Abdelaziz, F.~Mohamed, S.~Mekhamer, and M.~Badr, ``Distribution system
  reconfiguration using a modified tabu search algorithm,'' \emph{Electric
  Power Systems Research}, vol.~80, no.~8, pp. 943--953, 2010.

\bibitem{lavorato2012imposing}
M.~Lavorato, J.~F. Franco, M.~J. Rider, and R.~Romero, ``Imposing radiality
  constraints in distribution system optimization problems,'' \emph{IEEE
  Transactions on Power Systems}, vol.~27, no.~1, pp. 172--180, 2012.

\bibitem{el2016optimal}
S.~A. El~Batawy and W.~G. Morsi, ``Optimal secondary distribution system design
  considering rooftop solar photovoltaics,'' \emph{IEEE Transactions on
  Sustainable Energy}, vol.~7, no.~4, pp. 1662--1671, 2016.

\bibitem{li2014distribution}
J.~Li, X.-Y. Ma, C.-C. Liu, and K.~P. Schneider, ``Distribution system
  restoration with microgrids using spanning tree search,'' \emph{IEEE
  Transactions on Power Systems}, vol.~29, no.~6, pp. 3021--3029, 2014.

\bibitem{delage2010distributionally}
E.~Delage and Y.~Ye, ``Distributionally robust optimization under moment
  uncertainty with application to data-driven problems,'' \emph{Operations
  research}, vol.~58, no.~3, pp. 595--612, 2010.

\bibitem{xiong2017distributionally}
P.~Xiong, P.~Jirutitijaroen, and C.~Singh, ``A distributionally robust
  optimization model for unit commitment considering uncertain wind power
  generation,'' \emph{IEEE Transactions on Power Systems}, vol.~32, no.~1, pp.
  39--49, 2017.

\bibitem{bian2015distributionally}
Q.~Bian, H.~Xin, Z.~Wang, D.~Gan, and K.~P. Wong, ``Distributionally robust
  solution to the reserve scheduling problem with partial information of wind
  power,'' \emph{IEEE Transactions on Power Systems}, vol.~30, no.~5, pp.
  2822--2823, 2015.

\bibitem{qiu2015distributionally}
F.~Qiu and J.~Wang, ``Distributionally robust congestion management with
  dynamic line ratings,'' \emph{IEEE Transactions on Power Systems}, vol.~30,
  no.~4, pp. 2198--2199, 2015.

\bibitem{bagheri2017data}
A.~Bagheri, J.~Wang, and C.~Zhao, ``Data-driven stochastic transmission
  expansion planning,'' \emph{IEEE Transactions on Power Systems}, vol.~32,
  no.~5, pp. 3461--3470, 2017.

\bibitem{magnanti1995optimal}
T.~L. Magnanti and L.~A. Wolsey, ``Optimal trees,'' \emph{Handbooks in
  operations research and management science}, vol.~7, pp. 503--615, 1995.

\bibitem{baran1989network}
M.~E. Baran and F.~F. Wu, ``Network reconfiguration in distribution systems for
  loss reduction and load balancing,'' \emph{IEEE Transactions on Power
  delivery}, vol.~4, no.~2, pp. 1401--1407, 1989.

\bibitem{gao2017resilience}
H.~Gao, Y.~Chen, S.~Mei, S.~Huang, and Y.~Xu, ``Resilience-oriented
  pre-hurricane resource allocation in distribution systems considering
  electric buses,'' \emph{Proceedings of the IEEE}, 2017.

\bibitem{arab2015stochastic}
A.~Arab, A.~Khodaei, S.~K. Khator, K.~Ding, V.~A. Emesih, and Z.~Han,
  ``Stochastic pre-hurricane restoration planning for electric power systems
  infrastructure,'' \emph{IEEE Transactions on Smart Grid}, vol.~6, no.~2, pp.
  1046--1054, 2015.

\bibitem{yamangil2014designing}
E.~Yamangil, R.~Bent, and S.~Backhaus, ``Designing resilient electrical
  distribution grids,'' \emph{arXiv preprint arXiv:1409.4477}, 2014.

\bibitem{sa2002reliability}
Y.~Sa, ``Reliability analysis of electric distribution lines,'' Ph.D.
  dissertation, McGill University, 2002.

\bibitem{zhao2017distributionally}
C.~Zhao and R.~Jiang, ``Distributionally robust contingency-constrained unit
  commitment,'' \emph{IEEE Transactions on Power Systems}, vol.~33, no.~1, pp.
  94--102, 2018.

\bibitem{zeng2013solving}
B.~Zeng and L.~Zhao, ``Solving two-stage robust optimization problems using a
  column-and-constraint generation method,'' \emph{Operations Research
  Letters}, vol.~41, no.~5, pp. 457--461, 2013.

\bibitem{mccormick1976computability}
G.~P. McCormick, ``Computability of global solutions to factorable nonconvex
  programs: Part {I}--convex underestimating problems,'' \emph{Mathematical
  programming}, vol.~10, no.~1, pp. 147--175, 1976.

\end{thebibliography}


\end{document}